\documentclass[3p,12pt]{elsarticle}

\usepackage{amssymb}
\usepackage{amsthm}
\usepackage{amsmath}
\usepackage{empheq}
\usepackage{graphicx}
\usepackage{booktabs}
\usepackage{lineno}
\usepackage{subcaption}
\usepackage{algorithm}
\usepackage{multicol}
\usepackage{color}
\usepackage{enumitem}
\usepackage{lipsum}
\usepackage[utf8]{inputenc}

{\makeatletter\gdef\reallynopagebreak{\par\nopagebreak\@nobreaktrue}}

\makeatletter
\newcommand*{\NoBreakPar}{\vspace{\baselineskip}\par\nobreak\@afterheading}
\makeatother

\graphicspath{{allpic/}}
\newcommand\mynobreakpar{\par\nobreak\@afterheading}

\begin{document}
	\begin{frontmatter}
		\title{An extra-component method for evaluating fast matrix-vector multiplication with special functions}
		\author{Andrew V. Terekhov}
		\ead{andrew.terekhov@mail.ru}
		\address{Institute of Computational Mathematics and Mathematical Geophysics, 630090, Novosibirsk, Russia}
		\address{Novosibirsk State Technical University, 630073, Novosibirsk, Russia}
		\begin{abstract}
			In calculating integral or discrete transforms, use has been made of fast algorithms for multiplying vectors by matrices whose elements are specified as values of special (Chebyshev, Legendre, Laguerre, etc.) functions.  The currently available fast algorithms are several orders of magnitude less efficient than the fast Fourier transform. To achieve higher efficiency, a convenient general approach for calculating matrix-vector products for some class of problems is proposed. A series of fast simple-structure algorithms developed under this approach can be efficiently implemented with software based on modern microprocessors. The method has a pre-computation complexity of $O(N^2 \log N)$ and an execution complexity of $O(N \log N)$.  The results of computational experiments with the algorithms show that these procedures can decrease the calculation time by several orders of magnitude compared with a conventional direct method of matrix-vector multiplication.
		\end{abstract}
	\end{frontmatter}

\section{Introduction}
The discrete Fourier transform (DFT) has become a very popular method of numerical analysis due to the invention of \textcolor{black}{the fast Fourier transform algorithm (FFT) \cite{Cooley1965}.} The method decreases computational costs from $O(N^2)$ to $O(N\log N)$ in the calculation of matrix-vector products of the form
\begin{equation}
\begin{array}{l}
\mathbf{Y}=\mathcal{F}\mathbf{X} , \quad \mathcal{F}\in \mathbb{C}^{ N\times N}, \quad \mathbf{X},\mathbf{Y}\in \mathbb{C}^N, \\\\
\mathcal{F}:=\frac{1}{\sqrt{N}}\left(\exp{\left(-\frac{2\pi\mathrm{i}jk}{N}\right)}\right)_{j,k=0}^{N-1},
\end{array}
\label{fourier}
\end{equation}
where $\mathcal{F}$  is a Fourier transform matrix and $\mathrm{i}=\sqrt{-1}$. This algorithm has become a breakthrough in the development of methods of mathematical simulation and digital signal processing.

\textcolor{black}{Algorithms known as non-uniform fast Fourier transforms (NUFFTs), which allow fast calculation of transform (\ref{fourier}) for both unevenly specified samples and uneven sets of frequencies, have also been widely used in numerical analysis \cite{Dutt1993,Greengard2004,Plonka2018}. Most NUFFT methods are based on a procedure known as "gridding"{}~\cite{Jackson1991,Boyd1992,Dutt1993,Greengard2004,Beylkin1995,Liu1998,Fessler2003}. However, there are alternative approaches: Taylor expansion~\cite{Anderson1996}, fractional Fourier transform~\cite{Bailey1991}, low-rank approximation~\cite{RuizAntoln2018}, or the Butterfly algorithm~\cite{ONeil2010}.}
This naturally brings up the question of whether efficient procedures can be created not only for trigonometric functions, but also for classical orthogonal polynomials, such as Chebyshev, Legendre, Gegenbauer, Jacobi, Laguerre, and Hermite ones \cite{NIST:DLMF}. Let us consider some of the existing approaches to constructing such fast algorithms.

The economical method for the discrete Chebyshev transform, which can be calculated using the fast discrete cosine transform (DCT) based on the FFT algorithm, is one of the simplest methods based on classical orthogonal polynomials \cite{Ahmed1974,WenHsiungChen1977,Plonka2018,Driscoll2014}. The fast Chebyshev transform can be used to construct an efficient algorithm for Legendre polynomials \cite{Alpert1991,Hale2014} and, in a more general case, for Gegenbauer and Jacobi polynomials \cite{JieShen,Slevinsky2017,Bremer2019,Micheli2013}. Various approaches to calculating the expansion coefficients with a changed basis are studied in \cite{Townsend2016,Bostan2010,Tygert2010,Potts1998,Keiner2009}. The problem of implementing the fast Hermite transform is considered in \cite{Leibon2008}, where a fast algorithm with a computational complexity of $O(N\log^2N)$ operations is proposed.
Some more general algorithms are based on the idea of a preliminary compression  of the original matrix, for instance, by using the wavelet transform \cite{Beylkin1991}, the local cosine transform \cite{Aharoni1993,Matviyenko1996,Mohlenkamp1999}, or the Butterfly algorithm \cite{Michielssen1996,ONeil2010,Yin2019,Wedi2013}. \textcolor{black}{As a rule, a large amount of arithmetic operation is needed} for the preliminary compression of the matrix, but the method allows decreasing the total calculation time for multiple calculations of matrix-vector products. In the present paper, a simple and efficient method of matrix compression based on the standard FFT procedure will be considered for some class of algorithms.

{\color{black}
The remainder of this paper is organized as follows: In Subsection~2.1, a procedure of compression of transformation matrices is considered. In Subsection~2.2, an extra-component method is proposed for calculating matrix-vector products for trigonometric functions. In Section~3, to implement the fast Jacobi transform, a block version of the extra-component method is discussed. A fast algorithm to calculate the forward Laguerre transform is investigated in Subsection~4.1 and to calculate the backward one, in Subsection~4.2. The results of computational experiments are described in Section~5. Section~6 concludes the paper. In Appendix A, an efficient modification of the extra-component method is given for trigonometric functions.
}

\section{Extra-component algorithm}
\subsection{Non-uniform fast trigonometric transform}
Consider Chebyshev polynomials  $T_n:\mathbb{R}\rightarrow \mathbb{R}$ of the first kind defined on the interval  $\Omega=[-1,1]$
\begin{equation*}
\label{cheb_first_kind}
T_n(x)=\cos(n\arccos(x)),\; x\in \Omega,\ n \in \mathbb{Z}^+,
\end{equation*}
which form an orthogonal basis in  $L_{2,\omega}(\Omega)$
$$
\int_{-1}^{1}T_n(x)T_m(x)\omega(x)dx=\frac{\pi}{c_n}\delta_{mn},
$$
where $\omega(x)=\left(1-x^2\right)^{-1/2}$, $c_0=1$,  $c_n=2$ for $n\geq 1$, and $\delta_{mn}$ is the Kronecker delta.
Let a function  $f:\Omega\rightarrow \mathbb{R}$  be given, and let there exist the integral
$$
\|f\|^2_{L_{2,\omega}(\Omega)}=\int_{-1}^{1}\omega(x)|f(x)|^2dx.
$$
Then there exists a representation of the form
\begin{equation}
\label{sum_cheb1}
f(x)=\sum_{k=0}^{\infty}\hat{f}_k T_k(x), \quad x\in \Omega,
\end{equation}
\begin{equation*}
\label{int_cheb1}
\hat{f}_k=\frac{c_k}{\pi}\int_{-1}^{1}f(x)T_k(x)\frac{dx}{\sqrt{1-x^2}}=\frac{c_k}{\pi}\int_{0}^{\pi}\cos(k\theta)f(\cos(\theta))d\theta.
\end{equation*}
Assume that the values of $\hat{f}_k$ are known, and consider a problem of calculating a partial sum for (\ref{sum_cheb1}) which is written in the form of a matrix-vector product $\mathbf{F}=A\hat{\mathbf{F}}$, where
\begin{equation}
\quad A:=\left(\cos(m\arccos(x_n))\right)_{n,m=0}^{N,M}\in \mathbb{R}^{(N+1)\times( M+1)},
\label{cheb_matrix}
\end{equation}
$$
\mathbf{F}=\left(f(x_n)\right)_{n=0}^{N}, \quad \hat{\mathbf{F}}=\left(\hat{f}_m\right)_{m=0}^{M}.
$$
For a Chebyshev set of nodes, $ x_n=\cos(n\pi/N)$, and $M=N$ the matrix-vector multiplication (\ref{cheb_matrix}) can be calculated in $O(N\log N)$ arithmetic operations using the fast DCT \cite{Boyd2001}. \textcolor{black}{In what follows, a new algorithm will be developed to calculate $A\hat{\mathbf{F}}$ and $A^{\mathrm{T}}\mathbf{F}$ with $O(N\log N+M\varrho(\varepsilon))$ and $O(M\log M+N\varrho(\varepsilon))$ arithmetic operations, respectively, where $\varrho(\varepsilon) \leq 25$ is a function of the required accuracy.
}

\begin{figure}[!htb]
\centering
\includegraphics[width=\textwidth]{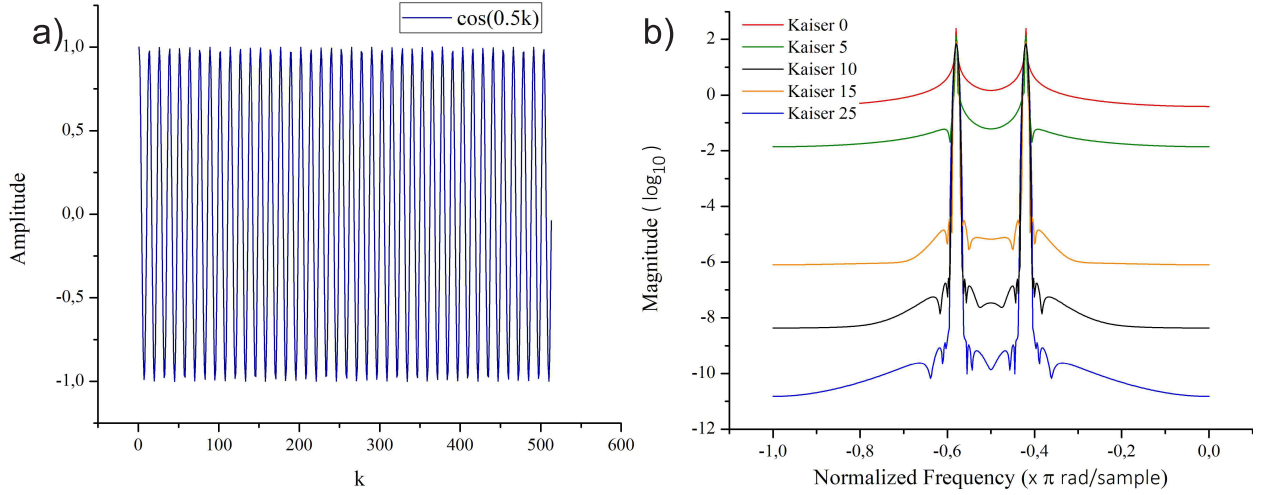}
\caption{a) row $A_{1\cdot}$ with points connected by lines; b) absolute values of the coefficients of the corresponding Fourier series for various window functions
}
\label{cheb_row}
\end{figure}

\begin{figure}
  \centering
  \includegraphics[width=\textwidth]{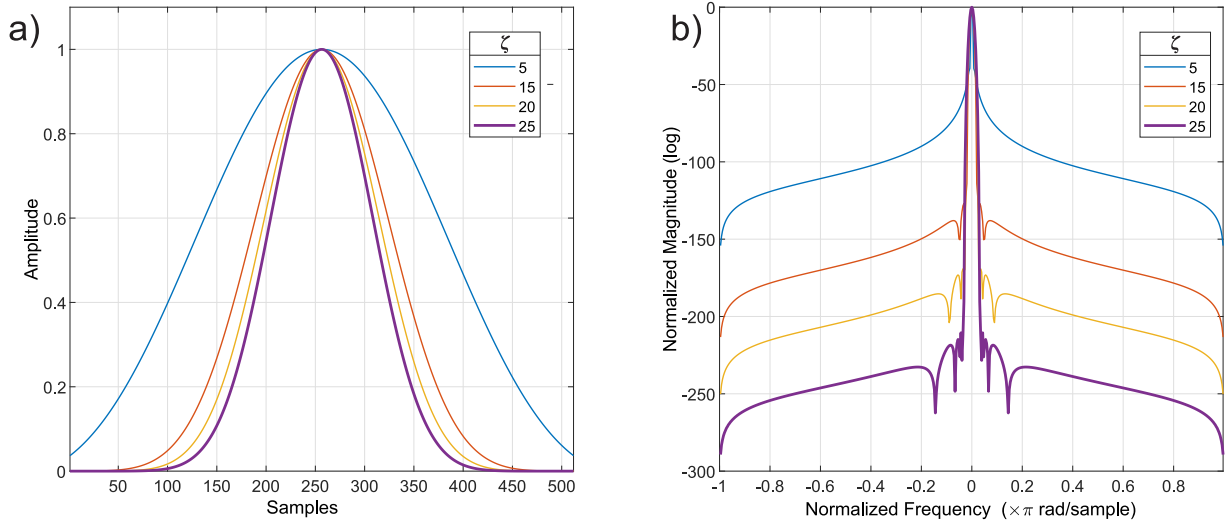}
  \caption{Kaiser window function for several values of the parameter $\zeta$ in a) time and b) frequency domains}
  \label{kaiser_win}
\end{figure}

Consider a row of the matrix $A$:
\begin{equation}
\label{cheb_row1}
A_{1\cdot}=\left[1,\cos(\theta_0),\cos(2\theta_0),...,\cos(M\theta_0)\right],
\end{equation}
where $\theta_0={1}/{2}$.
Fig.~\ref{cheb_row}a presents an example of the sequence $A_{1\cdot}$, and the absolute values of the corresponding Fourier components, $\mathcal{F}A_{1\cdot}^\mathrm{T}$, are shown in Fig.~\ref{cheb_row}b (see~curve "Kaiser 0"). Due to the sharp changes in the function at the boundaries, the absolute values of all coefficients of the Fourier series are nonzero, that is, the spectrum is not localized or, at least, finite. \textcolor{black}{In this case, the so-called frequency leakage effect is observed. \cite{Harris1978,Prabhu2018}.}
To eliminate this undesirable effect, one can use, for instance, the Kaiser window \cite{Prabhu2018}:
\begin{equation*}
\label{kaiser_function}
w^{\zeta,N}_n=\frac{I_0\left( \zeta \sqrt{1-\left(\frac{2n}{N}-1\right)^2}\right)}{I_0(\zeta)}, \quad 0 \leq n \leq N.
\end{equation*}
Here $I_0$ is the zeroth-order modified Bessel function of the first kind, and $\zeta \geq 0$ determines the shape of the Kaiser window (Fig.~\ref{kaiser_win}). 

\begin{figure}[!htb]
\centering
\includegraphics[width=1\textwidth]{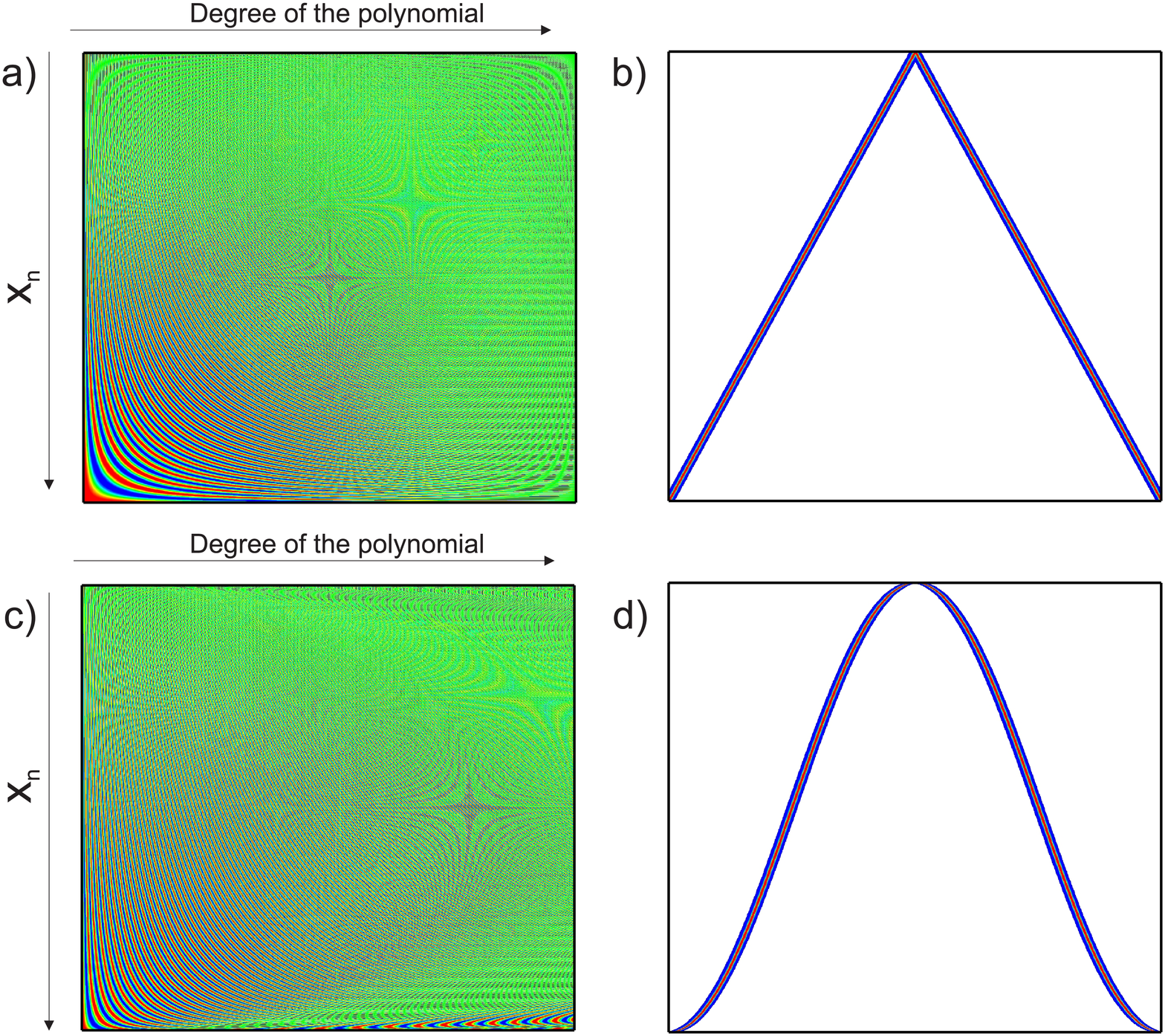}
\caption{$A\in \mathbb{R}^{1024\times 1024}$ for a) Chebyshev nodes and c) equispaced nodes. Absolute values of elements of the compressed matrix $AW_{1024}\mathcal{F}\in \mathbb{C}^{1024\times 1024}$ for b) Chebyshev nodes and d) equispaced nodes}
\label{pic:cheb_matrix22}
\end{figure}

Now consider a transform $\ddot{A}_{1\cdot}=\mathcal{F}W_MA_{1\cdot}^\mathrm{T}$ for the row (\ref{cheb_row1}), where $W_M$ is a diagonal matrix defined as
$$W_{M}=\mathrm{diag}\left\{w^{\zeta,M}_0,w^{\zeta,M}_1 ,...,w^{\zeta,M}_M\right\}\in \mathbb{R}^{(M+1)\times( M+1)}.$$
The absolute values of the spectrum components are presented in Fig.~\ref{cheb_row}b in a logarithmic scale. As the parameter $\zeta$ increases, the spectrum for (\ref{cheb_row1}) becomes to a large extent localized so that for some $\varepsilon>0$ the matrix elements satisfying the condition $|\tilde{a}_{1j}|<\varepsilon\max_j|\tilde{a}_{1j}|$ may be considered zero. A similar transform applied to all rows of the matrix $A$ provides a compressed matrix
\begin{equation}
\ddot{A}=AW_M\mathcal{F},
\label{A_comp}
\end{equation}  in which a significant number of elements can be ignored, since they are relatively small (see~Fig.~\ref{pic:cheb_matrix22}).

Taking into account the property of orthogonality,  $\mathcal{F}\mathcal{F}^{*}=\mathcal{F}^{*}\mathcal{F}=I$, where $I$ is the unit matrix, we can write the product $\mathbf{F}=A\hat{\mathbf{F}}$ as
\begin{equation}\label{comp1}
\mathbf{F}=A\underbrace{W_M\mathcal{F}\mathcal{F}^*W_M^{-1}}_I\hat{\mathbf{F}}=\ddot{A}\mathcal{F}^*W_M^{-1}\hat{\mathbf{F}},
\end{equation}
and the product  $\hat{\mathbf{F}}=A^{\mathrm{T}}\mathbf{F}$  as
\begin{equation}\label{comp2}
\hat{\mathbf{F}}=\underbrace{W_N^{-1}\mathcal{F}^*\mathcal{F}W_N}_I A^{\mathrm{T}}\mathbf{F}=W_N^{-1}\mathcal{F}^*\ddot{A}^{\mathrm{T}}\mathbf{F}.
\end{equation}
Thus, the preliminary calculation of the matrices  $\ddot{A}$ or $\ddot{A}^{\mathrm{T}}$, which have compact sparsity patterns, provides an efficient calculation of the desired matrix-vector products. For each vector to be multiplied, only one Fourier transform and one multiplication by a compressed matrix ($\ddot{A}$ or $\ddot{A}^{\mathrm{T}}$) are required. The multiplication by the diagonal matrix $W_M$ increases the total number of calculations only slightly.
\subsection{Reducing the computational errors}
\label{sec:reduce}
The efficient procedures for multiplying a vector by the matrix (\ref{cheb_matrix})  have been fully described, except for the calculation errors in formulas (\ref{comp1}) and (\ref{comp2}). The elements of the diagonal matrix $W_{N}^{-1}$ may become very large. Therefore, a matrix-vector product for this matrix may increase the errors caused by ignoring relatively small elements of the compressed matrices $\ddot{A}$  and $\ddot{A}^{\mathrm{T}}$.

Let us consider a modification of the calculation formulas (\ref{comp1}) to exclude the increase in the errors. For this an augmented matrix $A_e$ is obtained from the matrix $A$ by adding some new columns:
\begin{equation*}
A_e:=\left(\cos(m\theta_n)\right)_{n=0,m=-s}^{N,M+s},
\label{extd_f_e}
\end{equation*}
where $s\geq 0$ is a parameter of the number of columns added on the left and right sides.
Since the extra elements of the matrix $A_e$ are defined by formula (\ref{cheb_matrix}), the sparsity pattern of the matrix $\ddot{A}_e$ remains compact. Let us also form a new vector $\hat{\mathbf{F}}^{[1]}$ supplemented with zeros:
\begin{equation}\label{Fe1}
\hat{\mathbf{F}}^{[1]}=\left(\begin{array}{c}
\underbrace{
0 ,
...,
0}_{extra},
\underbrace{
  \hat{f}_0 ,
  \hat{f}_1 ,
  ... ,
  \hat{f}_M}_{\hat{\mathbf{F}}},
  \underbrace{
0 ,
  ...,
0}_{extra}
\end{array}
\right)^{\mathrm{T}},
\end{equation}
to be consistent with the dimensions of the augmented matrix $A_e$. Thus, when using formula (\ref{comp1}) the first and last $s$ of the diagonal elements of the matrix $W_{M+2s}^{-1}$ are multiplied by the fictitious zeros of the vector $\hat{\mathbf{F}}^{[1]}$. This allows controlling the accuracy of calculations by choosing the proper value of the parameter $s$.

To calculate the product $A^{\mathrm{T}}\mathbf{F}$  by formula  (\ref{comp2}) with accuracy control some extra rows are added:
\begin{equation*}
\hat{\mathbf{F}}^{[2]}=A_e^{\mathrm{T}}\mathbf{F}=\left(\cos(m\theta_n)\right)_{m=-s,n=0}^{M+s,N}\mathbf{F}.
\label{extd_f_e2}
\end{equation*}
Once the vector   $\hat{\mathbf{F}}^{[2]}$ is calculated, the first and last $s$ components are discarded:
\begin{equation}\label{Fe2}
\hat{\mathbf{F}}^{[2]}=\left(\begin{array}{c}
\underbrace{
\hat{{f}}_{-s} ,
...,
\hat{{f}}_{-1}}_{extra},
\underbrace{
  \hat{{f}}_0 ,
  \hat{{f}}_1 ,
  ... ,
  \hat{{f}}_M}_{\hat{{\mathbf{F}}}},\underbrace{
  \hat{{f}}_{M+1} ,
  ...,
  \hat{{f}}_{M+s}}_{extra}
\end{array}
\right)^{\mathrm{T}}.
\end{equation}

As mentioned above, the parameter $\zeta$ determines the form of the Kaiser window (Fig.~\ref{kaiser_win}) and controls the accuracy and number of the calculations. The degree of matrix compression depends on the parameter $\zeta$, since the multiplications of the rows of the matrix $A$ by the Kaiser window can limit the sharp changes between the first and last elements of the rows and, thus, localize the spectrum. 
The parameter $\zeta=\zeta(\varepsilon_1)$  can be found by Newton's method \cite{GALANTAI200025} for solving the equation $1/I_0(\zeta)=\varepsilon_1$  with an initial approximation $\zeta_0=50$.  \textcolor{black}{
The parameter $s=s(\zeta,\varepsilon_2)$ is chosen as small as possible so that \mbox{$w_s^{\zeta,M+2s}<\varepsilon_2$}. This can be achieved by a sequential exhaustive search over $s=1,2,3...$ by setting $\varepsilon_2>\varepsilon_1$ and precalculating $\zeta(\varepsilon_1)$. The constant $\varepsilon_1$ determines the elements of the compressed matrix, which should be regarded as zero since their values are small. The constant $\varepsilon_2$  determines the number of the extra components which are added to limit the growth of error caused by the multiplication by the matrix $W_{M+2s}^{-1}$. Thus, the parameters $\zeta(\varepsilon_1)$  and $s(\varepsilon_2,\zeta)$ are calculated using a small number of arithmetic operations to calculate matrix-vector products with an accuracy proportional to $\varepsilon_1/\varepsilon_2$}.

Let us now formulate some fast algorithms of multiplying  matrices  $A$ and $A^{\mathrm{T}}$ of the form (\ref{cheb_matrix})  by a vector as follows:

\begin{algorithm}[H]
\caption{Multiplication $\mathbf{F}=A\hat{\mathbf{F}}$}
\label{alg1}
{
\normalsize
\begin{enumerate}
\item \textit{Precomputation stage  }

  \subitem 1.1 For a given  $\varepsilon_1<\varepsilon_2$ calculate $\zeta(\varepsilon_1)$ and $s(\varepsilon_2,\zeta)$.
  \subitem 1.2 Set the extended matrix  $A_e:=\left(\cos(m\theta_n)\right)_{n=0,m=-s}^{N,M+s}$.
  \subitem 1.3 Calculate the compressed matrix $\ddot{A}_e=A_eW_{M+2s}\mathcal{F}$ and store the elements with absolute values greater than~$\varepsilon_1\|\ddot{A}_e\|_{max}$.
    \item \textit{Computation stage  }
   \subitem 2.1 Set $\hat{\mathbf{F}}^{[1]}$  according to (\ref{Fe1}).
   \subitem 2.2 Calculate $\mathbf{F}=\ddot{A}_e^{}\mathcal{F}^*W_{M+2s}^{-1}\mathbf{\hat{F}}^{[1]}$.
\end{enumerate}
}
\end{algorithm}

The following algorithm is used to calculate the forward Chebyshev transform:
\begin{algorithm}[H]
\caption{ Multiplication  $\hat{\mathbf{F}}=A^\mathrm{T}\mathbf{F}$}\label{alg2}
{
\normalsize
\begin{enumerate}
\item \textit{Precomputation stage}
 \subitem  1.1 For a given  $\varepsilon_1< \varepsilon_2$ calculate  $\zeta(\varepsilon_1)$ and $s(\varepsilon_2,\zeta)$.
  \subitem 1.2 Set the extended matrix  $A^{\mathrm{T}}_e:=\left(\left(\cos(m\theta_n)\right)_{m=-s,n=0}^{M+s,N}\right)^{\mathrm{T}}$.
  \subitem 1.3 Calculate the compressed matrix  $\ddot{A}_e^{\mathrm{T}}=\mathcal{F}W_{M+2s}A_e^{\mathrm{T}}$ and store the elements with absolute values greater than~$\varepsilon_1\|\ddot{A}_e^{\mathrm{T}}\|_{max}$.
      \item \textit{Computation stage }
  \subitem 2.1 Calculate   $\hat{\mathbf{F}}^{[2]}=W_{M+2s}^{-1}\mathcal{F}^*\ddot{A}_e^{\mathrm{T}}\mathbf{F}$.
   \subitem 2.2 Discard the first and last  $s$ elements of the vector   $\hat{\mathbf{F}}^{[2]}$ from (\ref{Fe2}) to form the required vector $\hat{\mathbf{F}}$.
\end{enumerate}
}
\end{algorithm}
Algorithms~1 and~2 can be used for multiple matrix-vector products not only for matrix (\ref{cheb_matrix}), but also for fast multiplication of matrices \textcolor{black}{that consist of the elements $\sin(m\theta_n)$ and $\exp(\mathrm{i}m\theta_n)$}. \textcolor{black}{The precomputation stage for the extra-component method requires about $O(NM\log M)$ arithmetic operations, and the computation stage, $O(M\log M+\varrho(\varepsilon)N)$ operations. However, for trigonometric functions the cost of preprocessing can be considerably decreased (see Appendix A).}

\textcolor{black}{To calculate the matrix-vector products considered above, NUFFT algorithms \cite{Plonka2018}  based on the "gridding"{}  procedure can be used. Although both the extra-component algorithm and the NUFFT-type methods use the FFT procedure and the Kaiser function \cite{Fessler2003}, the approach proposed in this paper is somewhat different.}

\textcolor{black}{To implement the NUFFT based on "gridding" the sequence of samples must be $2\pi$-periodic, while a special modification is used for an arbitrary length interval (see, for instance, formulas (108)-(110) in \cite{Dutt1993}). In constructing the extra-component method, the property of $2\pi$-periodicity of trigonometric functions is not taken into account, and no additional constraints are imposed.}

\textcolor{black}{
The Kaiser function in the extra-component method and the "gridding"{} procedure are used for various purposes. The "gridding"{} procedure uses the Kaiser (or Gaussian) function for spreading the approximated function values from a non-uniform grid to an oversampled uniform grid. In contrast to this, the extra-component method uses the Kaiser function to eliminate the frequency leakage effect in the compression step of the transformation matrix. As a result, the compressed matrix $\ddot{A}$ will have a compact sparsity pattern. From the point of view of the "gridding"{} procedure the compressed matrix columns (see Fig.~\ref{pic:cheb_matrix22}) for Algorithm~2 can be considered as a set of individual spreading wavelets for each sample, by means of which the transition to an equispaced DFT takes place. However, this interpretation is only a particular case of the extra-component method. Actually, for matrices of the form
}

\begin{figure}[!htb]
\centering
\includegraphics[width=1\textwidth]{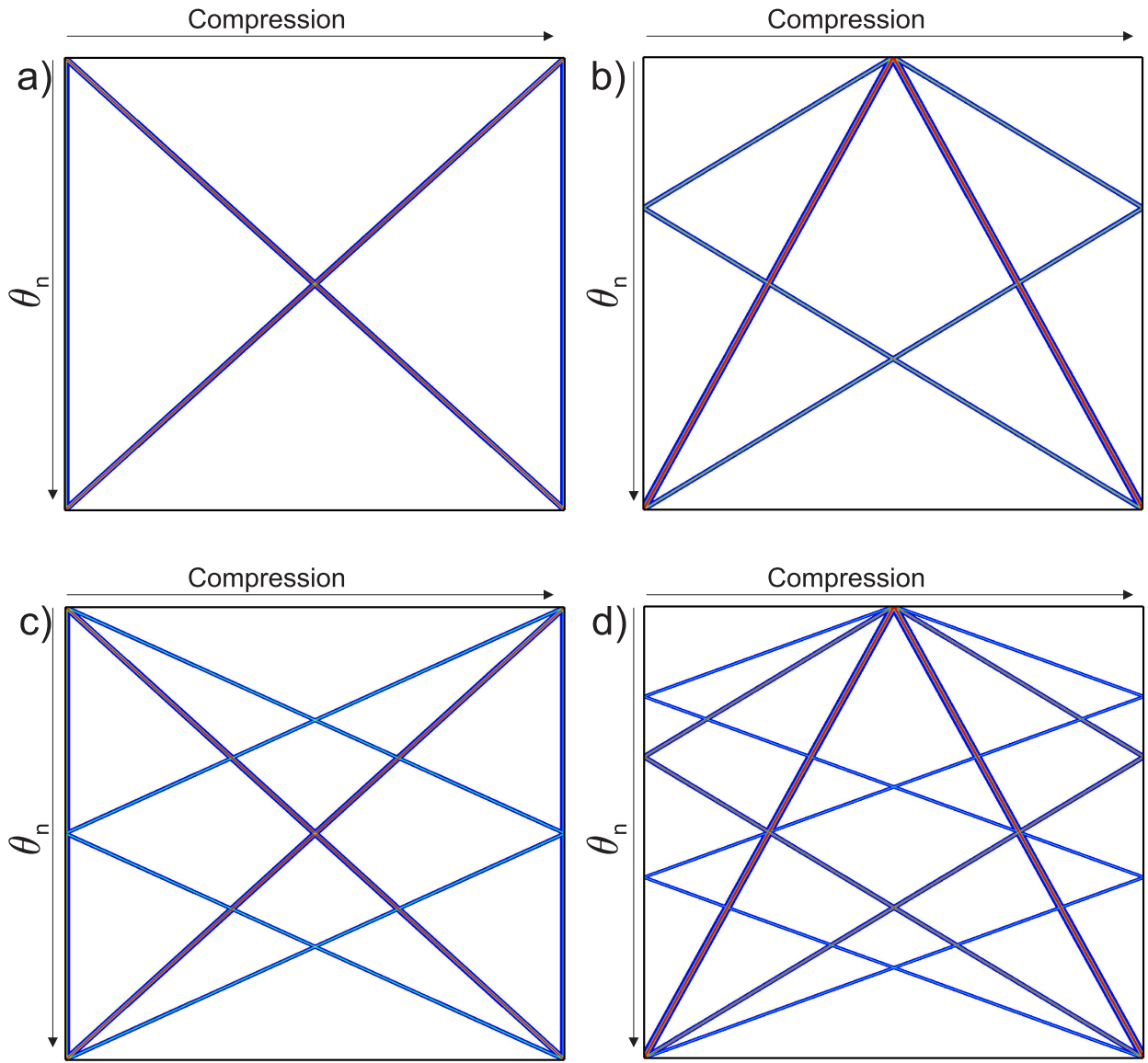}
\caption{\color{black}Absolute values of elements of the compressed matrix \mbox{$AW_{1024}\mathcal{F}\in \mathbb{C}^{1024\times 1024}$} for matrices of the form (\ref{cosd2}): a)~$\cos^2(m\theta_n)$, b)~$\cos^3(m\theta_n)$, c)~$\cos^4(m\theta_n)$, d)~$\cos^5(m\theta_n)$ and $\theta_n=\pi n/1024$ }
\label{pic:cheb_matrix_cosn}
\end{figure}

\textcolor{black}{
\begin{equation}\label{cosd2}
A:=\left(\cos^s(m\theta_n)\right)_{n=0,m=0}^{N,M},\  s=2,3,4,5....,
\end{equation}
the corresponding compressed matrices have compact sparsity patterns (Fig.~\ref{pic:cheb_matrix_cosn}). This makes it possible to use Algorithms~1 and~2 without additional modifications for the multiplication of matrices of the form  (\ref{cosd2}). "Gridding"{}-type methods are less universal and, therefore, more general algorithms (for instance, the Butterfly algorithm) are used in this case~\cite{ONeil2010}. However, as will be shown by computational experiments, the transformation matrix is compressed more efficiently by using the extra-component method. This results in a smaller value of the proportionality constant in the estimated algorithmic complexity of the method.}

\section{\textcolor{black}{Block extra-component algorithm}}
 Let a function $f(x)$ be given on the interval $\Omega=[-1,1]$ and square integrable with the Jacobi weight:
$$
\|f\|^2_{L_{2,\omega^{(\alpha,\beta)}(\Omega)}}=\int_{-1}^{1}\omega^{(\alpha,\beta)}(x)|f(x)|^2dx, \quad \quad \omega^{(\alpha,\beta)}(x)=(1-x)^{\alpha}(1+x)^{\beta},
$$
where $\alpha,\beta>-1$.
Then the function can be presented as a series
\begin{equation*}
\label{sum_jacobi}
f(x)=\sum_{m=0}^{\infty}\hat{f}^{(\alpha,\beta)}_m J_m^{(\alpha,\beta)}(x),\quad x\in \Omega,
\end{equation*} with expansion coefficients of the form
\begin{equation}
\begin{array}{c}
\label{int_jacobi}
  \displaystyle \hat{f}^{(\alpha,\beta)}_m=\frac{1}{\chi_m^{(\alpha,\beta)}}\int_{-1}^{1}\omega^{(\alpha,\beta)}(x)J_m^{(\alpha,\beta)}(x)f(x)dx,
\end{array}
\end{equation}
where $\chi_m^{(\alpha,\beta)}=\|J^{(\alpha,\beta)}_m\|_{L_{2,\omega^{(\alpha,\beta)}(\Omega)}}$.

\begin{figure}[!htb]
\centering
\includegraphics[width=\textwidth]{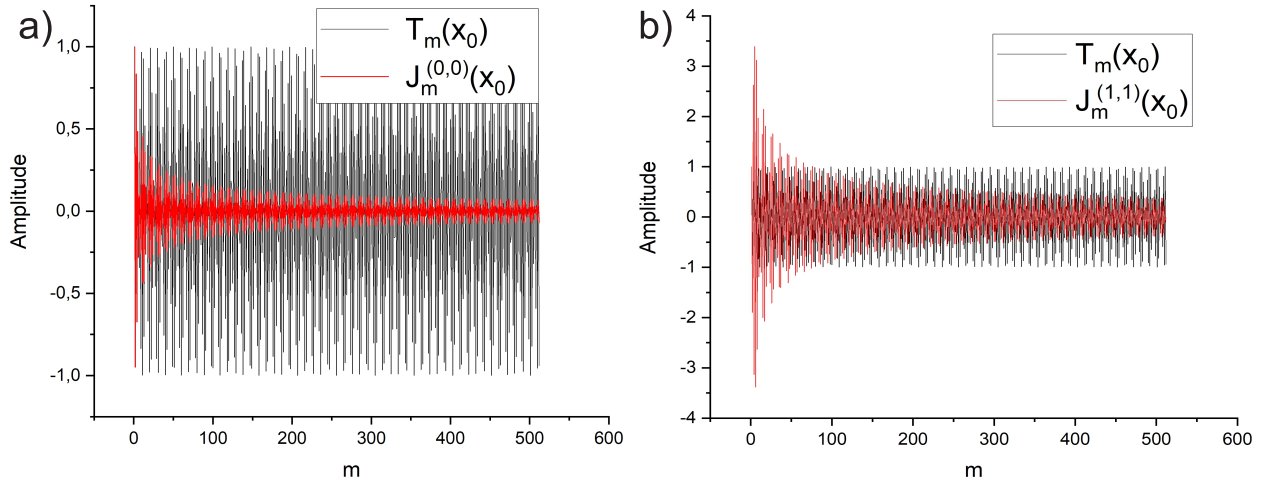}
\caption{a) Chebyshev and Legendre polynomials  and b) Chebyshev and Jacobi polynomials  of various orders at $x_0=1/2$}
\label{cheb_leg_1d}
\end{figure}

By specifying various values of the parameters $\alpha$ and $\beta$, one can obtain particular cases of Jacobi polynomials, namely: Chebyshev, Legendre, or Gegenbauer polynomials \cite{NIST:DLMF}. In what follows, a block version of the extra-component method will be proposed for calculating a matrix-vector product with the following matrix:
\begin{equation}
\quad B:=\left(J_m^{(\alpha,\beta)}(x_n)\right)_{n,m=0}^{N,M}\in \mathbb{R}^{(N+1)\times (M+1)}.
\label{Jac_matrix}
\end{equation}
Unfortunately, Algorithms~1 and~2 of the previous section cannot be directly applied to the matrix (\ref{Jac_matrix}) without modification, since, in contrast to the Chebyshev polynomials, the behavior of the Jacobi polynomials in a neighborhood $m=0$ changes considerably (Fig.~\ref{cheb_leg_1d}). It is not clear how the extra elements for the extended matrix $B_e$  should be specified at \mbox{$m<0$}  to make use of the procedures from Section~\ref{sec:reduce}.

\begin{figure}[!htb]
\centering
\includegraphics[width=\textwidth]{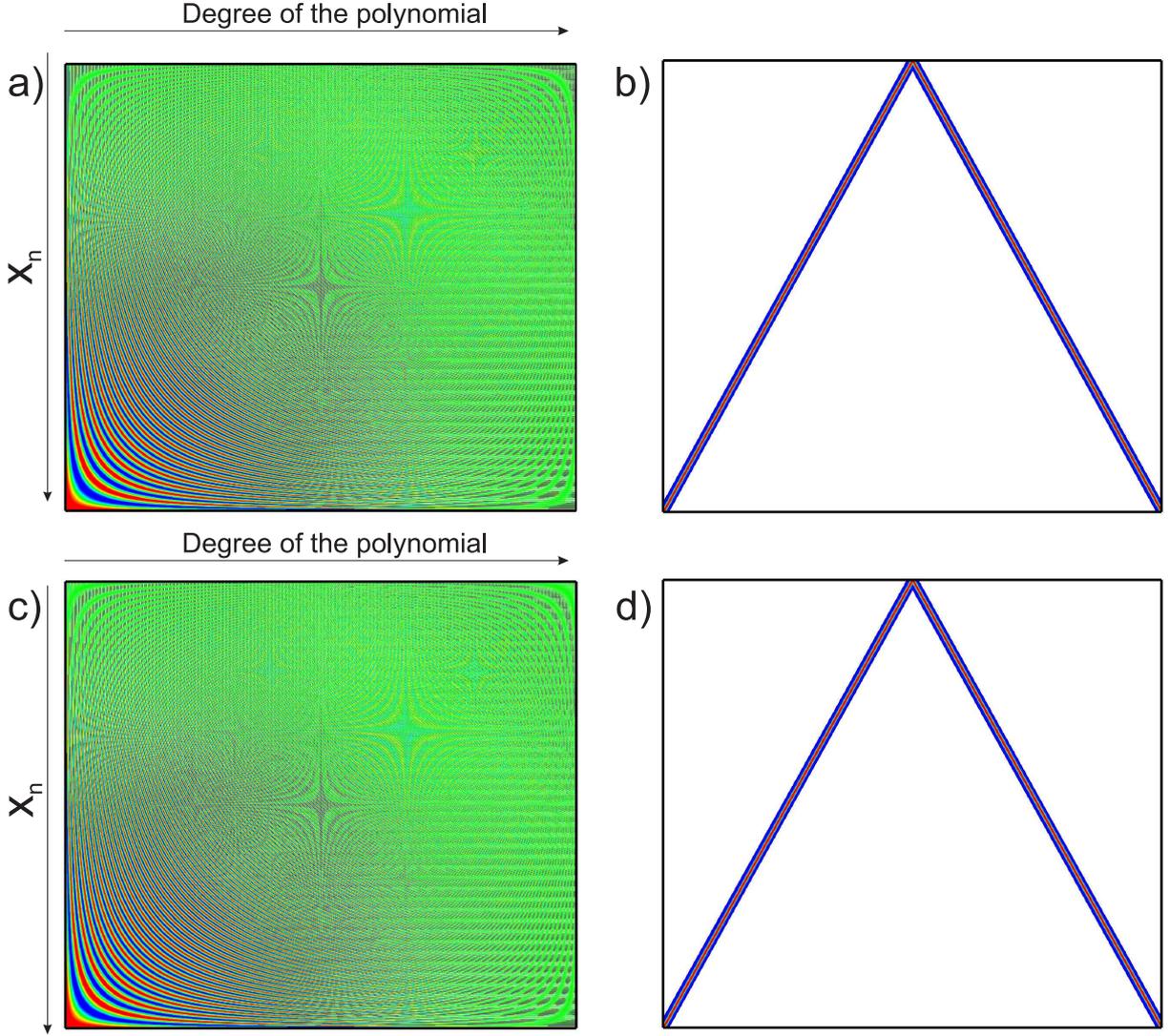}
\caption{ Matrix $B \in \mathbb{R}^{1024\times 1024}$ for a) $J^{(0,0)}_m$  and c) $J^{(1,1)}_m$. Absolute values of elements of the compressed matrix $BW_{1024}\mathcal{F}\in \mathbb{C}^{1024\times 1024}$ for polynomials b) $J^{(0,0)}_m$  and d) $J^{(1,1)}_m$ }
\label{pic:jac_comp1}
\end{figure}

\begin{figure}[!htb]
\centering
\subfloat{
\includegraphics[width={0.47\textwidth}]{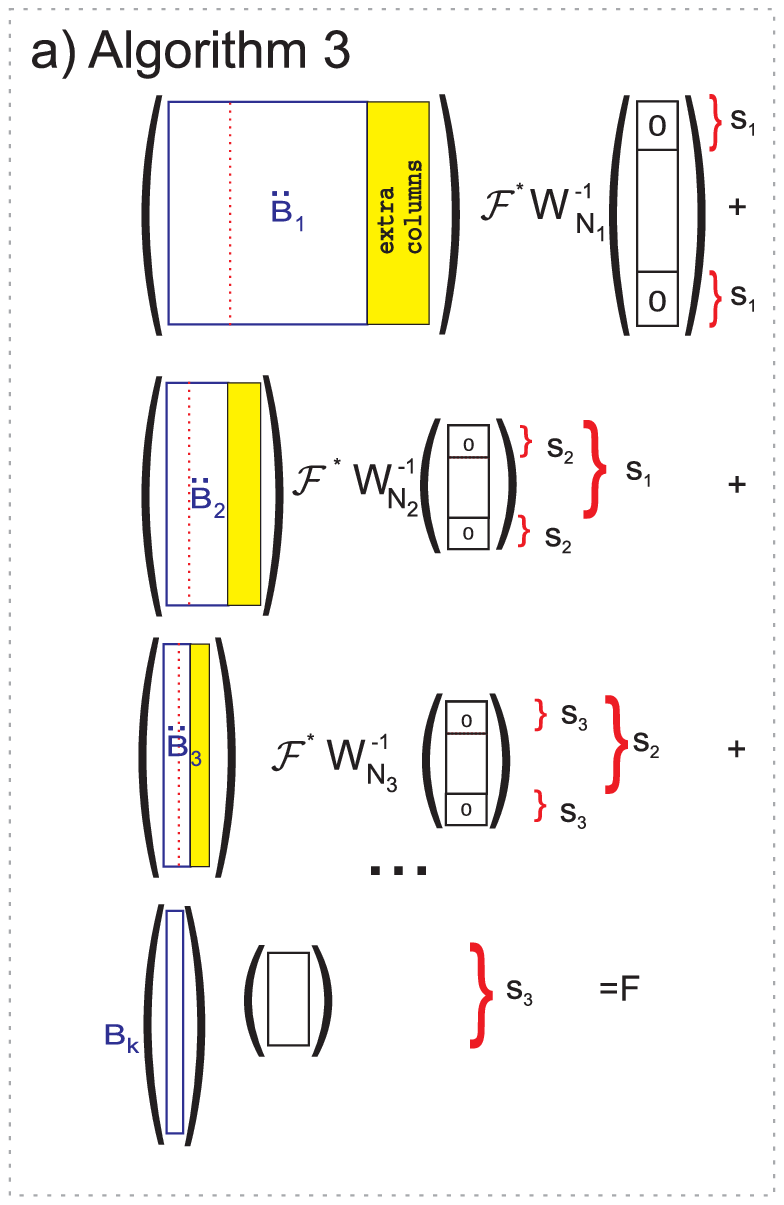}}\hfill
\subfloat{
\includegraphics[width={0.47\textwidth}]{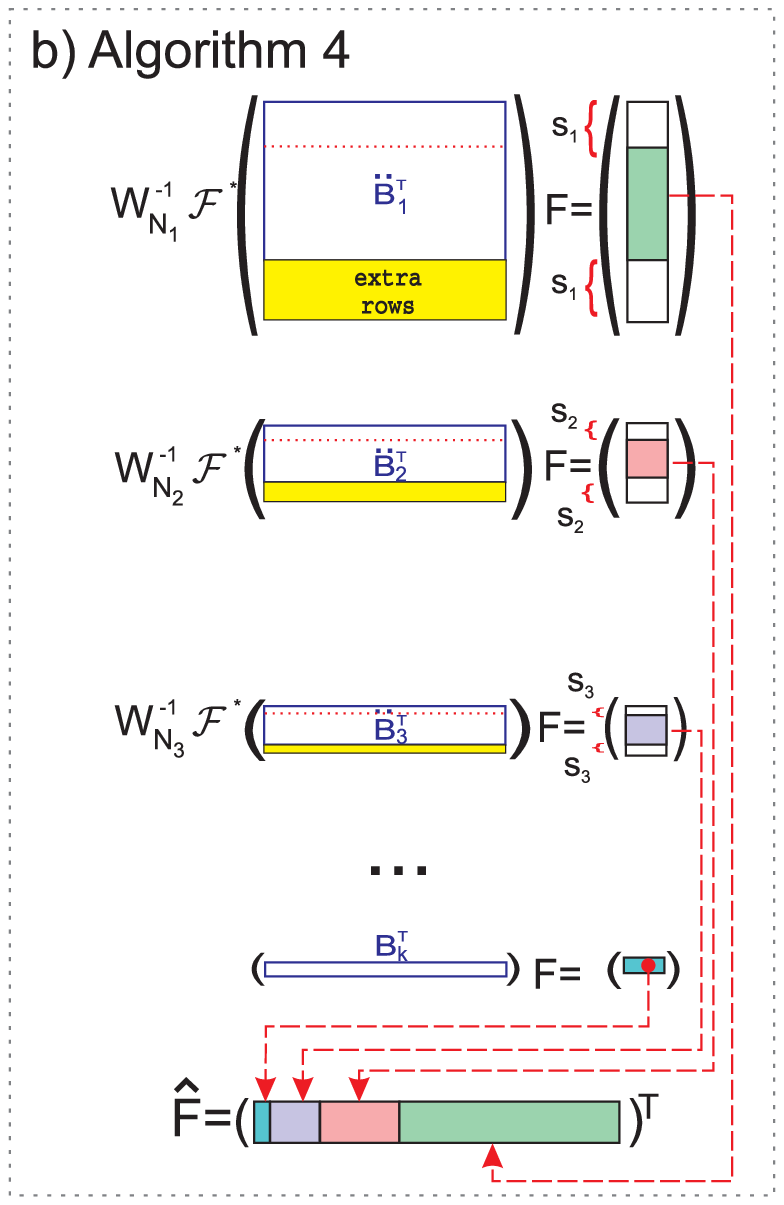}}
\hfill
\caption{Block extra-component method for calculating the product a)~$\mathbf{F}=B\hat{\mathbf{F}}$ and  b)~$\hat{\mathbf{F}}=B^{\mathrm{T}}\mathbf{F}$}
\label{pic:scheme200}
\end{figure}

Note, however, that if the transform (\ref{comp1})  is used for the matrix $B$, the sparsity pattern of the compressed matrix for Jacobi polynomials (Fig.~\ref{pic:jac_comp1}) does not differ from that for Chebyshev polynomials (Fig.~\ref{pic:cheb_matrix22}). Hence, a more general calculation process (Fig.~\ref{pic:scheme200}) can be proposed: \textcolor{black}{
to calculate $\mathbf{F}=B\hat{\mathbf{F}}$, matrix-vector products for each of the submatrices $B_i$ are calculated by Algorithm~\ref{alg1}, and for $\hat{\mathbf{F}}=B^\mathrm{T}\mathbf{F}$ by Algorithm~\ref{alg2}, respectively. Also, steps 1.2 for Algorithms~1 and~2 are modified as follows. Since no extra columns for the extended matrix $B_i$ can be added "on the left"{}, extra columns for the matrix $B_i$ are added only "on the right"{} by using $J_m^{(\alpha,\beta)}(x_n)$. The first $s_{i}$ columns of the initial matrix B are assumed to be extra columns of the matrix $B_i$.}
Therefore, the vector multiplied by the matrix is modified as follows:
\begin{equation*}
\label{Fe3}
\hat{\mathbf{F}}_i=\left(\begin{array}{c}
\underbrace{
0 ,
...,
0}_{extra},
  \hat{f}_{s_{i}+1} ,
  \hat{f}_{s_{i}+2} ,
  ... ,
  \hat{f}_M,
  \underbrace{
0 ,
  ...,
0}_{extra}
\end{array}
\right)^{\mathrm{T}}.
\end{equation*}
That is, the first $s_{i}$ components, corresponding to $ \hat{f}_{0},\hat{f}_{1},...,\hat{f}_{s_{i}}$  become zero, whereas in the extended vector  (\ref{Fe1})  zeros are added before these components. Hence, the thus-calculated $B_i\hat{\mathbf{F}}_i$  will contain no terms corresponding to the product of the first $s_{i}$ components of the vector $\hat{\mathbf{F}}$ and the first $s_{i}$ columns of the matrix $B$. These terms will be calculated at the next steps of the extra-component method, for which a submatrix consisting of the first $s_{i}$ columns of the initial matrix is formed. At the final step, the matrix-vector multiplication by the direct method is more efficient, since the size of the matrix at the last step is not large. The results of the calculations at the previous step are not used at the next step. Therefore, there arise no additional problems associated with stability. Of considerable interest is to study the degree of compression of all submatrices, since at $m=0$ the behavior of the Jacobi polynomials changes (Fig.~\ref{cheb_leg_1d}).~This question will be considered in the section devoted to computational experiments.

{ \color{black}
Fast algorithms to calculate matrix-vector products have been developed using precompression of the matrix by means of a wavelet (see \S8.1 in \cite{Mallat2008} ) or a block transformation (see \S8.3 in \cite{Mallat2008}). In the former case, the frequency axis is divided into approximation subintervals, while in the second case, the block transformation divides the time axis. An algorithm based on a wavelet transform for fast multiplication of vectors by matrices whose elements are non-oscillating functions is considered in \cite{Beylkin1991}, where a generalized Haar basis is taken as the wavelet. Although at first glance Fig.~\ref{pic:scheme200}b in the present paper and Fig.~1 in \cite{Beylkin1991} are similar, the extra-component method does not use multiresolution analysis based on the wavelet transform. With the extra-component method, a matrix is multiplied by dividing it into independent blocks whose size is determined by the parameters $\varepsilon_1$  and $\varepsilon_2$. To ensure the required accuracy, extra components (shown by yellow color in Fig.~\ref{pic:scheme200}) are added to each block.

Block discrete transforms are used in \cite{Mohlenkamp1999} to construct fast algorithms to expand a function as a series in associated Legendre functions (ALFs). Generally speaking, Algorithms~3 and~4 belong to this class of block methods. Their efficiency mostly depends on the choice of basis functions, transformation windows, and the strategy to determine the size and location of the approximation blocks. Algorithms~3 and~4 differ from the algorithms of \cite{Mohlenkamp1999} in that DFT is used instead of DCT and the Kaiser window is used to eliminate the frequency leakage effect. In the algorithms considered in \cite{Mohlenkamp1999} the strategy for selecting blocks for matrix compression is based on the dyadic decomposition method. In the extra-component method, the sizes of the nested computationally independent blocks are governed by the constants $\varepsilon_{1}$ and $\varepsilon_2$ determining the shape of the Kaiser window ($\zeta$) and the number of extra components ($s$) for each block. In \cite{Mohlenkamp1999}, orthogonal transformations are used to compress the matrices.
This method minimizes the computational errors in performing the inverse transformations. In the extra-component method, transform (\ref{A_comp}) is ill-conditioned. However, the accuracy of calculations when performing the inverse transformation in Algorithms~1--4 can be controlled by correctly setting the parameters $\varepsilon_{1}$ and $\varepsilon_2$. It will be shown by computational experiments that in the extra-component algorithm for the Legendre transform the compression of the transformation matrix is five to ten times better than in the method considered in \cite{Mohlenkamp1999}, and it also requires a smaller number of FFT computations.
}

{
\color{black}
Note that the approximation of local singularities based on global trigonometric functions may be inefficient; this is the case of ALFs of high orders, for which the range of function values may contain subdomains with an abrupt transition from a monotonic function behavior to an oscillating one. Algorithms~3 and~4 cannot be efficiently applied to high-order ALFs without modification, since the compressed matrix may have a large number of non-zero elements. In \cite{Mohlenkamp1999}, this problem is solved by an approximate analysis of the behavior of the ALF values to determine the optimal size of the matrix blocks for their subsequent compression. Unfortunately, the algorithm from \cite{Mohlenkamp1999} requires a larger computation time than the direct "on-the-fly"{} algorithms \cite{Schaeffer2013}. The extra-component method will be adapted and investigated for ALFs in a subsequent paper.
}

\section{\color{black} Application of the block extra-component algorithm to the Laguerre transform}
Consider the extra-component method when used for multiplying a vector by a matrix whose elements are defined as Laguerre functions:
\begin{equation*}
l_n(t)=\exp(-t/2)L_n(t),\quad t\geq 0,
\label{laguerre_function}
\end{equation*}
where $L_n(t)$ are Laguerre polynomials of degree  $n \in \mathbb{Z}^+$ \cite{NIST:DLMF}.
These functions form a complete orthonormal system of functions in  $ L_2[0,\infty)$:
\begin{equation*}
  \int_{0}^{\infty}l_m(t)l_n(t)dt=
  \left\{\begin{array}{ll}
  0,& m \neq n ,\\
  1,& m=n,
 \end{array}\right.
 \end{equation*}
such that for any function  $f(t) \in L_2[0,\infty)$ we have a representation in the form of a Laguerre series:
\begin{equation}
f(t)=\eta\sum_{m=0}^{\infty}\bar{f}_ml_m(\eta t),
\label{series_lag.sum}
\end{equation}
\begin{equation}
\bar{f}_m=\int_{0}^{\infty}f(t)l_m(\eta t) dt,
\label{series_lag.int}
\end{equation}
where $\eta>0$ is a parameter that controls the convergence rate of the series. The Laguerre transform is used in solving both forward and backward problems of mathematical modeling  \cite{Boyd2001,Terekhov2017,Terekhov2018,Mikhailenko1999,Abate1996}. Therefore, of interest is to reduce the number of operations when using formulas  (\ref{series_lag.sum}) and (\ref{series_lag.int}).

\subsection{Spectral-domain algorithm}
In the general case, to expand functions into a series in orthogonal polynomials it is necessary to calculate integrals of rapidly oscillating functions. For this, to ensure both stability and high accuracy of calculations, high-accuracy Gaussian quadratures for nonuniform grids can be used. However, in solving many problems the initial data are specified with a constant discretization step, which does not allow using high-accuracy Gaussian quadratures. A new approach for calculating the transform (\ref{series_lag.int}) has been developed in \cite{Terekhov2018a}. It is based on solving the transport equation by the classical method of separation of variables.

Let us formulate an auxiliary initial boundary value problem for the convection equation
\begin{equation*}
\begin{array}{ll}
\displaystyle \frac{\partial v}{\partial t }-\frac{\partial v}{\partial x }=0, \quad t>0,\quad  x\in[0,L],
\end{array}
\label{advection_eq}
\end{equation*}
 with conditions   $v(x,0)=f(x),\; v(0,t)=v(L,t)$. After the Laguerre transform with respect to time, this problem can be written as \cite{Mikh2003}
\begin{equation}
\label{advection_laguerre}
\left({\eta}/{2}-\partial_x\right)\bar{v}_m=-\Phi(\bar{v}_m),
\end{equation}
where  $\Phi(\bar{v}_m)=-f+\eta\sum_{j=0}^{m-1}\bar{v}_j.$
Since  $\Phi(\bar{v}_{m})=\eta \bar{v}_{m-1}+\Phi(\bar{v}_{m-1})$, from (\ref{advection_laguerre})  we obtain
\begin{subequations}
\label{advection_laguerre3}
  \begin{empheq}{align}
  \left({\eta}/{2}-\partial_x\right)\bar{v}_{0}-f=0,\\
  \left({\eta}/{2}-\partial_x\right)\bar{v}_{m}=\left(-{\eta}/{2}-\partial_x\right)\bar{v}_{m-1},\quad m=1,2,...
  \end{empheq}
\end{subequations}
Taking the Fourier transform with respect to the variable $x$, expressing the sought-for function in the spectrum, and going back to $x$, we obtain a solution to the problem (\ref{advection_laguerre3}) in the form \cite{Terekhov2018a}
\begin{equation}
\label{main_formula}
\displaystyle \bar{v}_{m}(x)=\sum_{j=0}^{\infty}\tilde{f}_j\frac{\left({-{\eta}/{2}-\mathrm{i}k_j}\right)^m}{\left({{\eta}/{2}-\mathrm{i}k_j}\right)^{m+1}}\exp\left(\mathrm{i}\frac{2\pi j x}{L}\right),
\end{equation}
where  $\tilde{f}_j$ are the Fourier series coefficients for the function  $f(x)$, which is given on the interval  $x\in[0,L]$ and  $k_j={2\pi} j/{L}$.
In accordance with the solution  (\ref{main_formula}), the function  $f(x)$, given as an initial condition, will "move"{} in the direction  $x_0=0$. Writing the solution at this point in terms of the Laguerre series coefficients, we obtain
\begin{equation}
\bar{v}_{m}=\bar{v}_{m}(0)=\sum_{j=0}^{\infty}\tilde{f}_j\frac{\left({-{\eta}/{2}-\mathrm{i}k_j}\right)^m}{\left({{\eta}/{2}-\mathrm{i}k_j}\right)^{m+1}}.
\label{main_formula2}
\end{equation}
This method of calculation, which is based on the classical method of separation of variables, adds fictitious periodicity of the form \mbox{$f(t)=f(t+bL)$}, where $b$ is an arbitrary nonnegative integer. To exclude the periodicity, two different approaches have been proposed \cite{Terekhov2018a}: one approach uses zero-padding, and the other one, the conjugation operation
\begin{equation}\label{conjg_operation}
\mathbb{Q}\left\{\bar{v}_j;\tau\right\}=\sum_{m=0}^{\infty}\left(\bar{v}_m-\bar{v}_{m-1}\right)l_{m+j}(\eta \tau),\ \quad j \in \mathbb{Z}^+, \ \bar{v}_{-1}\equiv0.
\end{equation}
It has been shown that if the conjugation operation, is applied twice $\mathbb{Q}^2\left\{\bar{v}_j;L\right\}\equiv \mathbb{Q}\left\{\mathbb{Q}\left\{\bar{v}_j;L\right\};L\right\}$, the fictitious periodicity on the approximation interval $[0,L]$ can be excluded. The conjugation (\ref{conjg_operation}) is a linear correlation of two sequences, and it can be calculated for a finite number of terms using the FFT algorithm  \cite{Nussbaumer1982}.

Assume that the function to be approximated is represented by a Fourier series with $N_x+1$ \textcolor{black}{ coefficients. Then the calculation} of the coefficients $\bar{v}_m$ by formula (\ref{main_formula2})  will require $O(N_xM)$ operations. It is easy to show that
 $$\displaystyle \frac{\left({-{\eta}/{2}-\mathrm{i}k_j}\right)^m}{\left({{\eta}/{2}-\mathrm{i}k_j}\right)^{m+1}}=\frac{\mathrm{i}\exp(-\mathrm{i}m\phi(k_j))}{{\eta}/{2}-\mathrm{i}k_j},$$
where $$\phi(k_j)=\arctan \left( -k_j,\eta/2 \right) -\arctan \left( -k_j,-\eta/2 \right),$$ and the function $\arctan(x, y)$ computes the principal value of the argument function of the complex number $x +\mathrm{i}y$. \textcolor{black}{Formula} (\ref{main_formula2})  can be written as a matrix-vector product \mbox{$\bar{\mathbf{V}}=C^{\mathrm{T}}\tilde{\mathbf{F}}$} , where
\begin{equation}\label{lag_matrix}
\begin{array}{l}\displaystyle
C:=\left(\frac{\mathrm{i}\exp(-\mathrm{i}m\phi(k_j))}{\eta/2-\mathrm{i}k_j}\right)_{j=0,m=0}^{N_x,M},\\\\ \displaystyle
 \mathbf{\bar{V}}=\left(\bar{v}_0,\bar{v}_1,...,\bar{v}_M\right)^{\mathrm{T}},\; \tilde{\mathbf{F}}=\left(\tilde{f}_0,\tilde{f}_1,...,\tilde{f}_{N_x}\right)^{\mathrm{T}}.
 \end{array}
\end{equation}
To multiply the matrices $C$ and $C^{\mathrm{T}}$, Algorithms~1 and~2 can be used without additional modifications, since the elements of the extended matrices for the function $\exp(-\mathrm{i}m\phi(k_j))$ are determined explicitly.
\begin{figure}[!htb]
\centering
\includegraphics[width=\textwidth]{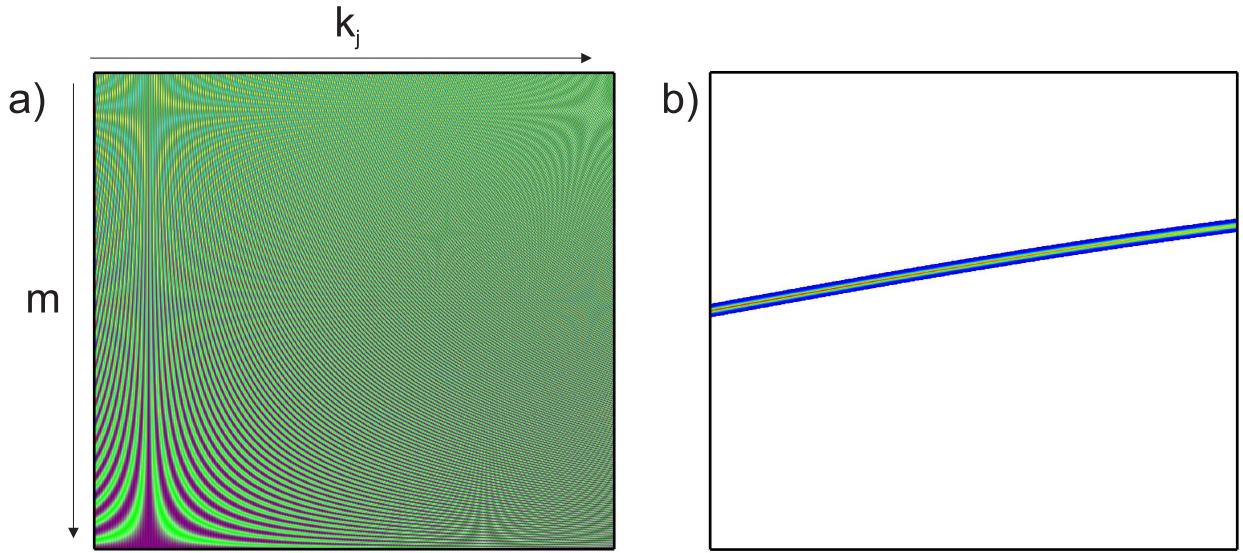}
\caption{a) Real part of the matrix $C^\mathrm{T} \in \mathbb{C}^{649\times 501}$ of the form (\ref{lag_matrix}) with parameters $\eta=2500$ and $k_j={2\pi j}/{T}$ for $T=4\ s$; b) Absolute values of the compressed matrix elements  $\mathcal{F}W_{649}C^\mathrm{T} \in \mathbb{C}^{649\times 501}$}
\label{lag_row}
\end{figure}
Fig.~\ref{lag_row} shows that once the procedure of compression is applied, the number of matrix elements that can be neglected increases considerably. The function $f(x)$ is assumed to be periodic with zero values, $f(0)=f(L)=0$, at the boundaries of the approximation interval $[0,L]$. A method for approximating functions of a more general form without this restriction is considered in \cite{Terekhov2018a}.

To calculate the inverse of the Laguerre transform (\ref{series_lag.sum}) using a spectral approach, the Laguerre series coefficients should be changed for the Fourier series coefficients for an equivalent approximation interval with no discontinuities of the function at the boundaries. In some cases this can be done on the basis of the operation (\ref{conjg_operation}). Nevertheless, an important question is whether fast Algorithms~1--4 can be used to calculate the sum (\ref{series_lag.sum}) directly in the time domain without using an auxiliary Fourier spectrum.
 \subsection{Time-domain algorithm}
\begin{figure}[!htb]
\centering
\includegraphics[width=\textwidth]{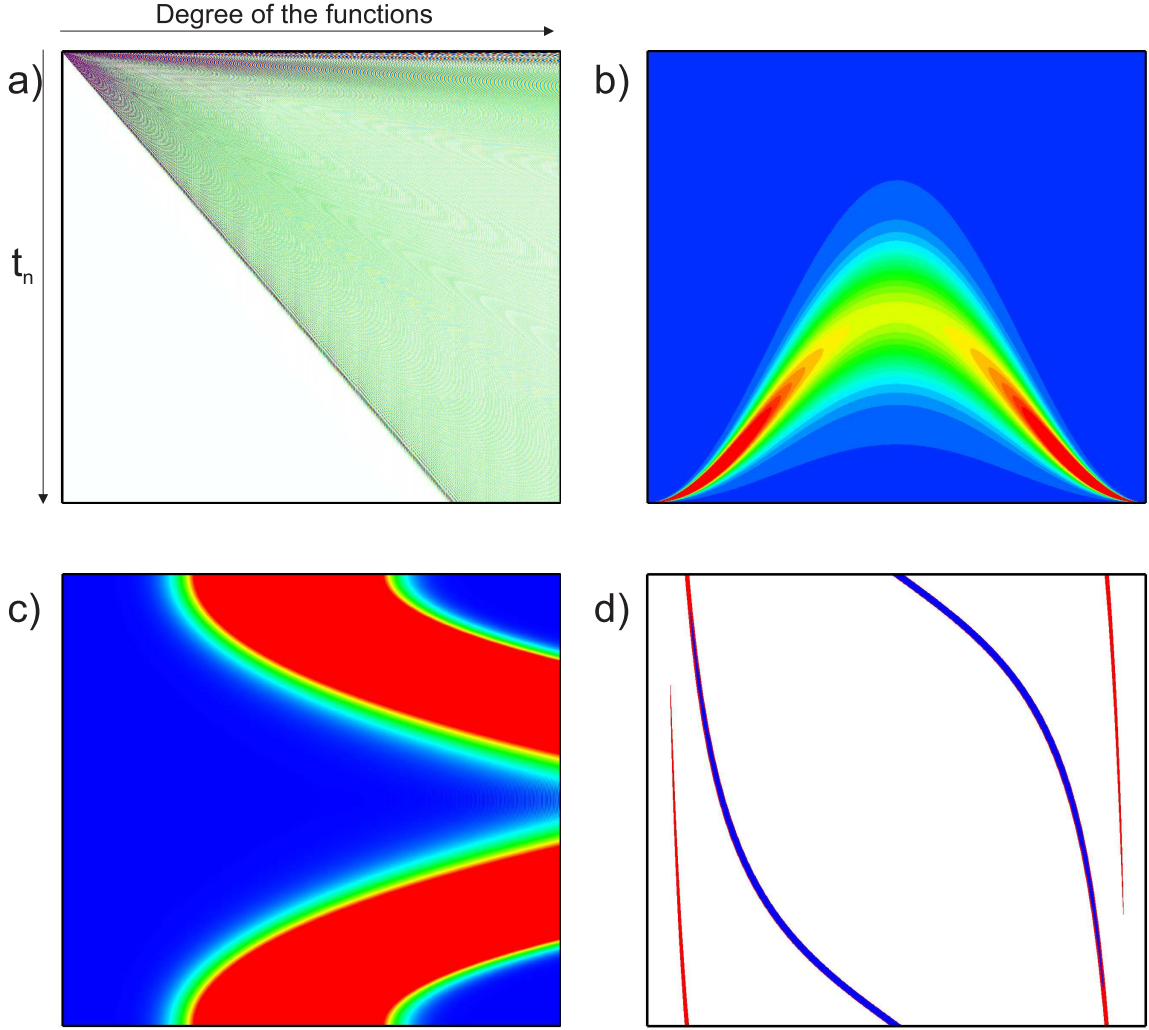}
\caption{a)~Laguerre transform matrix $D \in \mathbb{R}^{1000\times 1024}$ and compressed matrices b)~$D W_{1024}\mathcal{F}$, c)~$\mathcal{F}W_{1000}D$, and d)~$\mathcal{F}W_{1000}DW_{1024}\mathcal{F}$}
\label{pic:lag_matrix2}
\end{figure}

 For equispaced nodes, consider the transforms (\ref{comp1}) and (\ref{comp2}) for the matrix shown in Fig.~\ref{pic:lag_matrix2}a. Practical calculations show (Figs.~\ref{pic:lag_matrix2}b,c) that Algorithms~1--4 cannot be used efficiently, since the compression (\ref{comp1}) or (\ref{comp2}) does not decrease the number of nonzero matrix elements. However, if a transform of the form
 \begin{equation}\label{comp_2d}
 \ddot{D}=\mathcal{F}W_MDW_N\mathcal{F},
\end{equation}
is used, Fig.~\ref{pic:lag_matrix2}d shows that a considerable number of the compressed matrix elements can be neglected in comparison with the one-dimensional compression in Fig.~\ref{pic:lag_matrix2}b,c.

Calculations by formula (\ref{comp_2d}) can be made by using a combination of the algorithms in Fig.~\ref{pic:scheme200}a,b. For instance, at computation stage 2.1 of Algorithms~1 and~3, Algorithms~2 and~4 should be used for the multiplication by the matrix $\ddot{A}$  with some additional compression of the matrix $\ddot{A}$ rows. Another way is to use Algorithms~2 and~4 with the compression of the matrix rows and then its columns by using Algorithms~1 and~3. Thus, the matrix $\ddot{D}$  is a result of taking the two-dimensional Fourier transform of the matrix $W_MDW_N$. The two-dimensional compression is considered in more detail in Section \ref{par:lag}. To calculate the forward transform (\ref{series_lag.int}), it is not reasonable to abandon the spectral approach, since otherwise Gaussian quadratures with integration nodes defined on a nonuniform grid will have to be used, which may be inconsistent with the original data of the problem.

\section{Numerical experiments}
 Let us consider a series of computational experiments to estimate the efficiency of the above-proposed algorithms of fast matrix-vector multiplication for some special matrices.  All algorithms have been implemented as Fortran-2008 programs using BLAS and FFTW libraries of the Intel Math Kernel Library. \textcolor{black}{  The calculations have been made on the supercomputer of Novosibirsk State University. The supercomputer comprises Intel Xeon Gold~6248  $20$-core processors operating at $2.5$ GHz. Each computational node contains $4$ processors and 384 GB of RAM. }

\subsection{Trigonometric  transform }
\label{sect_num_tri}
Consider a matrix-vector multiplication of the form  (\ref{cheb_matrix})  to calculate the values of a Chebyshev series at points $x_n\in \Omega$. This class of matrices is widely used, since in the general case not only the Chebyshev transform, but also the Laguerre transform (as shown above), as well as many other applied problems of computational mathematics, can be reduced to a trigonometric basis. Since the elements of the matrix (\ref{cheb_matrix}) are real, the Fourier series coefficients are symmetrically conjugate. Hence, only a half of the compressed matrix needs to be calculated at the precomputation stage (Fig.~\ref{pic:cheb_matrix22}). The elements of the vector being multiplied will also be real. Therefore, the Fourier transform at the second stage of Algorithms~1--4  is performed to calculate half the spectrum of the vector for the dimensions to be consistent.

\begin{figure}[!htb]
\centering
\includegraphics[width={\textwidth}]{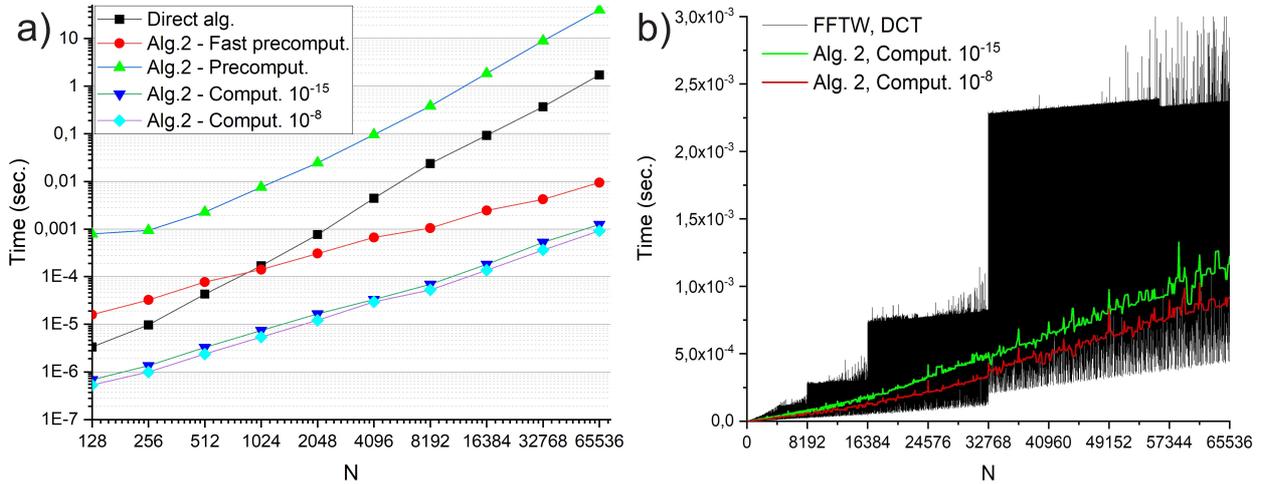}
\caption{a)~\textcolor{black}{Precomputation and computation times of non-uniform DCT for Alg.~2 and the direct method versus the transform size, b)~DCT computation time versus the transform size for Alg.~2 and fast DCT (fftw library)}}
\label{fig:cheb12}
\end{figure}

The results of the experiments are presented in Fig.~\ref{fig:cheb12}a, where the calculation time is defined as the averaged time of a thousand calculations. The matrix is applied to a vector containing random numbers uniformly distributed on the interval $[0, 1]$.  It is clear that the approach being proposed can significantly decrease the calculation time in a wide range of $N$-values. To achieve an accuracy of $10^{-8}$ the number of extra columns was specified as $s\approx N/4$, and to achieve an accuracy of $10^{-15}$ it was specified as $s\approx N/2$. In the first case the size of the original matrix increased by one and a half, and in the second case it increased by a factor of two, respectively. Regardless of the matrix order $N$, in the first case the number of stored diagonals of the compressed matrix (Fig.~\ref{pic:cheb_matrix22}) $bw=16$, and in the second case $bw=24$. \textcolor{black}{
Compared to the direct algorithm for calculating matrix-vector products, the extra-component method requires considerable precomputation. However, as shown in Appendix A, the precomputation time for trigonometric functions can be reduced from $O(N^2\log N)$ to $O(N \log N+Nq)$, $q<25$, which is in agreement with our calculation results (see "Fast precomput."{} vs. "Precomput." in Fig~\ref{fig:cheb12}a. )}

In the software implementation of the algorithms, to achieve maximum efficiency the parameter $s=s(\varepsilon_2,\zeta)$ (which affects not only the accuracy but also the efficiency of the extra-component method) must be chosen correctly. The computational complexity of the FFT algorithm is estimated at $O(N\log N)$. However, in practice (see~Fig.~\ref{fig:cheb12}b) the time of some fast algorithms (the DCT is implemented on the basis of the FFT) can differ considerably for neighboring values of $N$ and, hence, the above estimate becomes not valid. According to the computational experiments, the size of the extended matrix $A_e\in \mathbb{R}^{N,M+2s}$ should be such that the FFT could be made in the least time. For this, one could take a slightly larger parameter $s$ than required to achieve the required accuracy and, hence, select the locally optimal size of the Fourier transform. This approach calls for a preliminary assessment of the efficiency of the FFT procedure in a neighborhood of $M+2s$ values. \textcolor{black}{
As a result, the efficiency of the extra-component method seems to be rather high (see~Fig.\ref{fig:cheb12}b), and the DCT or DST (discrete sine transform) for nonoptimal transform sizes can be performed in less time.}

\subsection{Jacobi transform}
Using, as an example, the forward and backward Jacobi transforms, let us consider the efficiency of block Algorithms~3 and~4.
 { \color{black} With Gauss-Jacobi quadratures \cite{Gil2020} for discretization (\ref{int_jacobi}), the matrix $B \in \mathbb{R}^{N \times N}$ for the forward discrete Jacobi transform can be written as follows:
\begin{equation*}
\begin{array}{c}
B:=\left(\sqrt{ \frac{\omega^{(\alpha,\beta)}_n}{\chi^{(\alpha,\beta)}_m}}{J}_m^{(\alpha,\beta)}(x_n)f(x_n)\right)_{m,n=0}^{N-1},
\end{array}
\end{equation*}
where  $\omega^{(\alpha,\beta)}_n,\;x_n$ are weights and nodes of the quadrature formula.
The use of Gauss-Jacobi quadratures is not necessary in implementing the extra-component method, but in this case the matrix $B$  is orthogonal, that is, $B^{-1}=B^{\mathrm{T}}$.}

\begin{figure}[!htb]
\centering
\includegraphics[width=\textwidth]{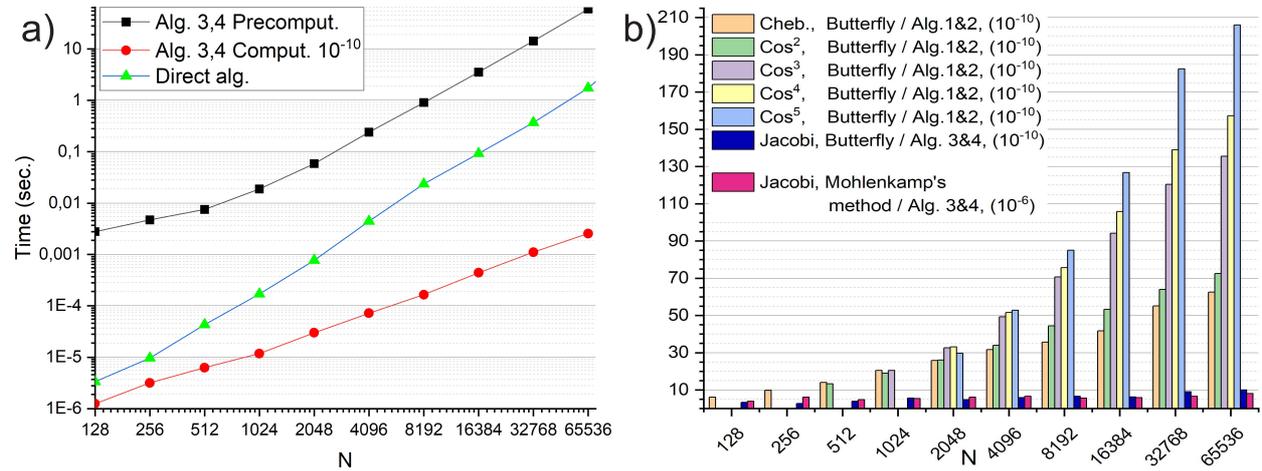}
\caption{a)~Precomputation and computation times versus the transform size  to perform forward (or backward) Legendre transform with an accuracy of $10^{-10}$, b)~\textcolor{black}{Ratio of the number of nonzero elements of the compressed matrix for the Butterfly and Mohlenkamp's algorithms to the number of nonzero elements of the compressed matrix for the extra-component method}}
\label{pic:jac1}
\end{figure}

\begin{figure}[!htb]
\centering
\includegraphics[width=\textwidth]{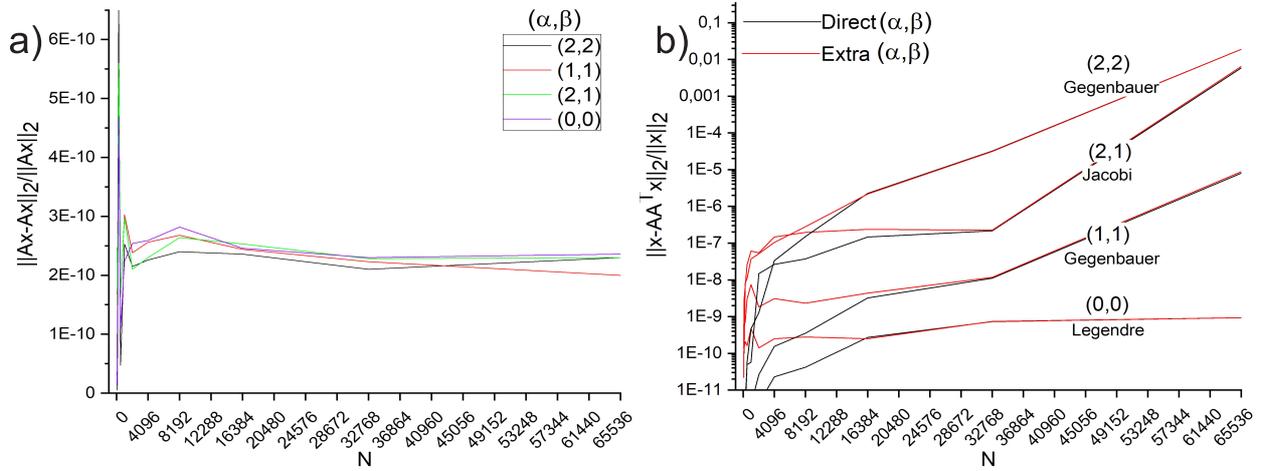}
\caption{Calculation accuracy of the Jacobi transform with parameters $(\alpha,\beta)$ versus the transform size for a) only forward transform and b) forward transform with subsequent backward transform}
\label{pic:test2}
\end{figure}

Figs.~\ref{pic:jac1} and~\ref{pic:test2}  show \textcolor{black}{the calculation time and accuracy} of the forward transform  and of a sequence of \textcolor{black}{forward and backward transforms versus} the parameters $\alpha$ and $\beta$. \textcolor{black}{In contrast to the direct method, the extra-component method (not taking into account the precomputation time) requires several orders of magnitude less computation time. Thus, if the transformation matrix is multiplied by different vectors multiple times, the extra-component algorithm is efficient. The calculations have shown that to provide an accuracy of $ 10^{-10}$ the compressed matrix must be stored with a bandwidth  $bw=20$ for $\varepsilon_1=10^{-10}$, $\varepsilon_2=10^{-3}$, whereas to achieve an accuracy of $ 10^{-6}$ we may take $bw=16$, $\varepsilon_1=10^{-6}$, $\varepsilon_2=10^{-4}$.}

With increasing $\alpha$ or $\beta$ the calculation accuracy of the forward transform does not change (Fig.~\ref{pic:test2}a). However, if the forward and backward transforms are used in sequence with increasing $\alpha$ or $\beta$, the error also begins to increase rapidly (Fig.~\ref{pic:test2}b), since $B B^{\mathrm{T}} \neq I$ due to the accumulation of errors when using the three-term recurrence relations that determine the values of the Jacobi polynomials.

\begin{table}[!htb]
  \centering
  \begin{tabular}{ccccccccc}
  \hline
  &\multicolumn{3}{c}{Alg.~3 ($10^{-6}$)} &\multicolumn{5}{c}{Alg.~3 ($10^{-10}$)}\\
  &\multicolumn{3}{c}{ $\varepsilon_1=10^{-6}$, $\varepsilon_2=10^{-4}$, bw$=16$} &\multicolumn{5}{c}{$\varepsilon_1=10^{-10}$, $\varepsilon_2=10^{-3}$, bw$=20$}\\
     \cmidrule(r){2-4}\cmidrule(r){5-9}
  N &$N_1$ & $N_2$ & $N_3$ & $N_1$ &$N_2$ &$ N_3$ &$N_4$ & $N_5$ \\
    \cmidrule(r){1-1} \cmidrule(r){2-4}\cmidrule(r){5-9}
 128 & 135 & - & - & 151 & 29 & - & - & - \\
  256 & 271 & - & - & 304 & 60& - & - & - \\
  512 & 539 & 33 & - & 600 & 160 & - & - & - \\
  1024 & 1078 & 60 & - & 1200 & 208 & 39 & - & - \\
  2048 & 2160 & 12 & - & 2340 & 336 & 51 & - & - \\
  4096 & 4312 & 230 & - & 4680 & 666 & 95 & - & - \\
  8192 & 8624 & 456 & 29 & 9360 & 1326 & 181 & 29 & - \\
  16384 & 17248 & 910 & 51 & 18720 & 2662 & 374 &58 & - \\
  32768 & 34496 & 1820 & 99 & 37440 & 5292 & 704 & 98 & - \\
  65536 & 68796 & 3430 & 179 & 74536 & 10192 & 1352 & 184 & 29 \\
  \hline
\end{tabular}
   \caption{\color{black}
   Order of submatrices $B_i \in \mathbb{R}^{N\times N_i}$ to be compressed to calculate the forward Legendre transform of dimension $N$ by Alg.~3 with an accuracy of $10^{-6}$ and $10^{-10}$.}
      \label{table_1}
\end{table}

\textcolor{black}{
In \cite{Slevinsky2017} and \cite{Bremer2019}, algorithms were proposed for calculating the Jacobi and Jacobi-Chebyshev transforms for $|\alpha|,|\beta| < 1/2$. Computational experiments have shown that Algorithms~3 and~4 can be used for a wider range of the parameters, $\alpha,\beta>-1$. Note that in the extra-component algorithm the costs of the precomputation step are $O(N^2\log N)$, which is several orders of magnitude greater than the costs of the computation step (Fig.~\ref{pic:jac1}). In comparison with the Butterfly  and Mohlenkamp's algorithms, the extra-component method requires a comparable amount of precomputation, but has a better compression of the transformation matrix (Fig.~\ref{pic:jac1}b). The transformation matrix in Algorithms~3 and~4 has smaller compression than that in Algorithms~1 and~2 due to the need to compress several submatrices (see Table~\ref{table_1}). For every $i$-th $N\times N_i$  submatrix the number of nonzero elements for the corresponding compressed matrix does not depend on $N_i$. It is of the order of $bw\cdot N$  for the Jacobi transform. The degree of compression of the entire transformation matrix determines to a great extent the efficiency of the computation step in the extra-component method. The FFT procedure requires no more than $10-20\%$ of the total computation time, and for each block the FFT dimension is $N_i$.}

\subsection{Laguerre transform}
Let us assess the efficiency of the above-proposed algorithms for multiplying a vector by a matrix defined by the series (\ref{main_formula2}). The software implementation of the extra-component method for this matrix does not differ from those for the Chebyshev transform or cosine transform discussed in Section~\ref{sect_num_tri}. Fig.~\ref{fig:cheb12}a shows that the extra-component method decreases the calculation time by several orders of magnitude in comparison to the direct method, and demonstrates high accuracy. Thus, if a function represented by a Laguerre series can be approximated by a rapidly converging Fourier series, the backward transform can be efficiently calculated using the spectral approach with Algorithm~1. With Algorithm~2, by solving the transport equation, the forward transform can be made, with subsequent exclusion of the fictitious periodicity \cite{Terekhov2018a}. 

\begin{figure}[!htb]
\centering
\includegraphics[width=\textwidth]{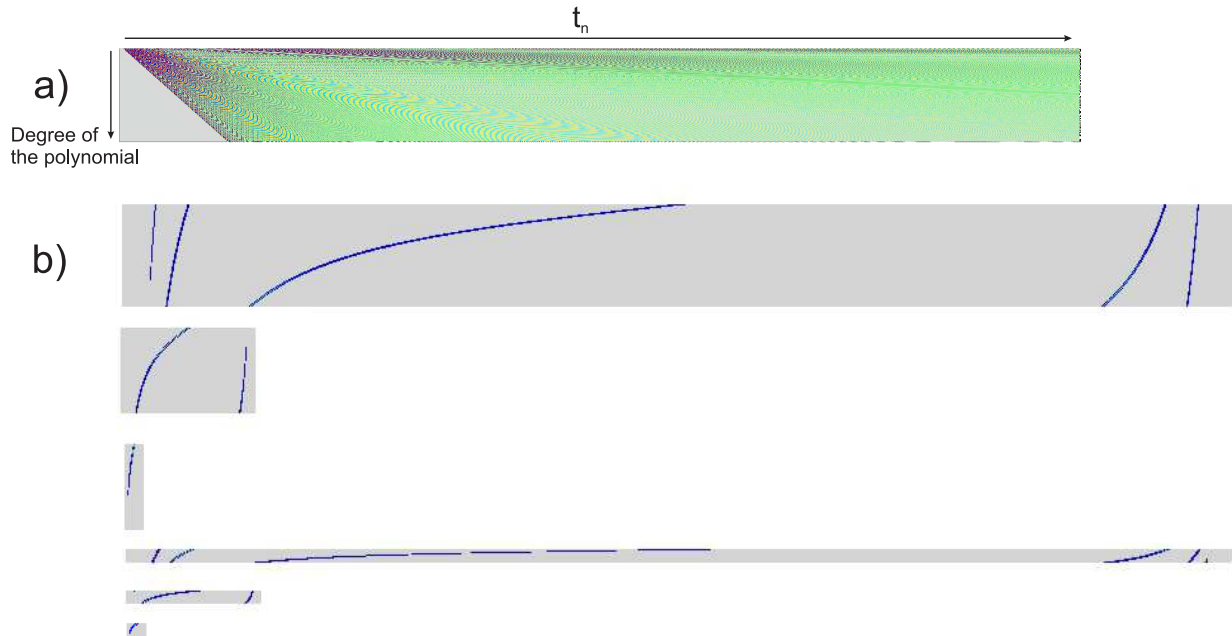}
\caption{ \color{black} Decomposition of a)~matrix of the backward Laguerre transform (\ref{series_lag.sum}) into b)~compressed submatrices in the 2D block version of the extra-component method}
\label{pic:test322}
\end{figure}

\begin{figure}[!htb]
\centering
\includegraphics[width=\textwidth]{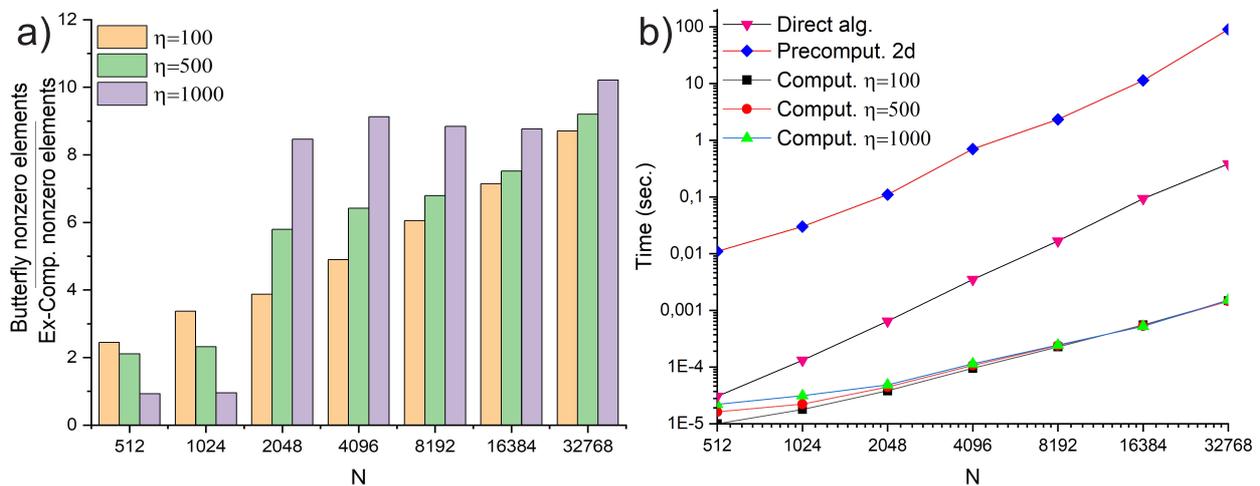}
\caption{a)~\textcolor{black}{
Ratio of the number of nonzero elements of the compressed matrix for the Butterfly algorithm to the number of nonzero elements of the compressed matrix for the extra-component method, b) precomputation and computation times for the extra-component and direct algorithms versus the matrix size to perform backward Laguerre transform with an accuracy of $10^{-10}$}}
\label{pic:test422}
\end{figure}

\label{par:lag}
Consider the  two-dimensional procedure of compression (\ref{comp_2d})  of the extra-component method for calculating the Laguerre series (\ref{series_lag.sum}) \textcolor{black}{
for $\eta=100,500 $ and $1000$. The order of the transformation matrix is set to be $N\times N$ and the Laguerre series values are calculated at points $t_i=12i/{(N-1)},\; i=1,2,...,N$.}
In contrast to the multiplication by the matrix  (\ref{main_formula2}), for the Laguerre transform inversion it will be necessary to multiply \textcolor{black}{by the submatrices shown in } Fig.~\ref{pic:test322}b. Without  the division by the Kaiser function, the final result could be obtained by multiplying only the first matrix by the vector. However, to exclude division by near-zero elements of the matrices $W_N$ and $W_M=W_N$ (see~Section~\ref{sec:reduce}), the multiplication will have to be performed \textcolor{black}{
performed by additional matrices of various orders, which (since their sparsity patterns are compact) can be efficiently multiplied by the vector}.  The transformation matrix can be compressed up to $3-4\%$ of the initial number of elements. \textcolor{black}{
As for Algorithms~3 and~4, the computation time of the FFT will be about an order of magnitude less than the time of multiplying the compressed submatrices by the subvectors. Therefore, the number of nonzero elements of the compressed matrix is one of the major factors affecting the efficiency of the method as a whole.  Fig.~\ref{pic:test422}a shows that the extra-component algorithm requires less operations for multiplying the compressed matrix than the Butterfly algorithm.}

The calculation times are given in Fig.~\ref{pic:test422}b, which shows that the two-dimensional algorithm is somewhat less efficient than the one-dimensional version of the extra-component method. This can be explained  by the fact that the calculations are performed using a  two-dimensional block scheme and also by the greater number of extra components added in comparison with the one-dimensional compression, since the errors of multiplication by the matrices $W_N^{-1}$ and $W_M^{-1}$ must be excluded. Nevertheless, the proposed two-dimensional method of matrix compression for calculating values of the series (\ref{series_lag.sum}) significantly reduces the calculation time in comparison with the direct multiplication algorithm without introducing any fictitious periodicity. Thus, although it is difficult to create a fast algorithm of the Laguerre transform that could compete with the direct method of matrix-vector multiplication \textcolor{black}{(especially for small transform size)}, some new ideas for solving this problem have been proposed in the present paper. This opens up additional possibilities for practical use of integral and discrete transforms in problems of numerical analysis.

\subsection{Non-oscillating functions}
{\color{black}
Let us consider matrices whose elements are given by non-oscillating functions with local singularities. In this case the extra-component method using a global trigonometric basis to compress the transformation matrix may be inefficient because the Fourier series may have slow convergence. Methods based on the discrete wavelet transform are more suitable for such matrices \cite{Beylkin1991}. Let us compare the degree of compression and the accuracy of matrix-vector multiplication by the extra-component method and the wavelet method for some test matrices \cite{Beylkin1991}.
\begin{equation*}
A^{[1]}_{i,j}=\left\{
\begin{array}{cc}
  \frac{1}{\pi} \Lambda(j), &  0= i  \leq j <N,\\\\
  \frac{2}{\pi} \Lambda(j-i) \Lambda(j+i), & 0< i  \leq j <N,\\\\
  0, &\text {otherwise},
\end{array}
\right.
\end{equation*}
where $ \Lambda(z)=\mathrm{\Gamma}(z+\frac{1}{2})/\mathrm{\Gamma}(z+1)$;

\begin{equation*}
A^{[2]}_{i,j}=\left\{
\begin{array}{cc}
   \displaystyle \frac{i\cos(\log i^2)-j\cos(\log j^2)}{(i-j)^2}, &  i \neq j,\\\\
  0, & i=j,
\end{array}
\right.
\end{equation*}

\begin{equation*}
A^{[3]}_{i,j}=\left\{
\begin{array}{cc}
  \displaystyle \frac{1}{i-j+\frac{1}{2}\cos(ij)} , &  i \neq j,\\\\
  0, & i=j.
\end{array}
\right.
\end{equation*}
To compress the matrices, the block version of the extra-component method based on transform (\ref{comp_2d}) was used. Table~\ref{table:wavelet} shows that for the matrices $A^{[1]}$  and $A^{[2]}$ both algorithms demonstrate a high degree of matrix compression. The maximum accuracy of the wavelet algorithm for the matrix $A^{[2]}$ does not exceed $10^{-2}$, whereas the accuracy of the extra-component method can be increased to $10^{-6}$ if a larger number of elements of the compressed matrix is used.

The extra-component method for the matrix $A^{[3]}$ did not allow compressing the transformation matrix to an acceptable level. This can be explained by the fact that the elements of the matrix $A^{[3]}$ are given by a function with a complex behavior of its singular values. In this case the Fourier series has a low convergence rate. With the wavelet transform such local singularities can be better taken into account.

\begin{table}[!htb]
\centering
\begin{tabular}{cccccc}
  \hline
  &\multicolumn{2}{c}{Matrix $A^{[1]}$}&\multicolumn{3}{c}{Matrix $A^{[2]}$}\\
  \cmidrule(r){2-3}      \cmidrule(r){4-6}
  N & Wavelet ($10^{-5}$) & Extra ($10^{-5}$)& Wavelet ($10^{-2}$) & Extra ($10^{-2}$) &   Extra ($10^{-6}$)   \\
      \cmidrule(r){1-1}   \cmidrule(r){2-2} \cmidrule(r){3-3}    \cmidrule(r){4-4}  \cmidrule(r){5-5}  \cmidrule(r){6-6}
  64 & 1.73 & 1.36 & 2.37 & 2.0 & 1.0 \\
  128 & 2.89 & 2.63 & 4.13 & 3.7 & 1.6  \\
  256 & 5.18 & 5.16 & 8.25 & 6.6 & 2.0  \\
  512 & 9.7 & 9.5 & 14.8 & 10.3 & 3.7  \\
  1024 & 18.6 & 19.6 & 33 & 23 & 6.5 \\
  \hline
\end{tabular}
\caption{\color{black} Accuracy of calculating matrix-vector multiplication and the ratio of the number of nonzero elements of the compressed matrix to the total number of elements $N^2$ for the wavelet-based algorithm and the extra-component method}
\label{table:wavelet}
\end{table}

}
\section{Summary and conclusions }
A new approach to constructing fast algorithms for the calculation of matrix-vector products related to integral and discrete transforms has been proposed.  The method developed has made it possible to decrease the calculation time by several orders of magnitude   for  various types of transforms, such as Chebyshev, Legendre, Gegenbauer, Jacobi  and  Laguerre.
In the precomputation step, which is performed only once, the initial matrix by using the fast Fourier transform is reduced to a matrix with a compact sparsity pattern. Then, in the computation stage, the matrix-vector product can be calculated rather quickly, since the near-zero elements of the compressed matrix are neglected.

{
\color{black}
With the precomputation step, it makes sense to use the extra-component method for multiple calculations of matrix-vector multiplication with the same matrix for different vectors. If the individual components of the Fourier spectrum for columns or rows of the original transformation matrix can be calculated efficiently in a number of operations that does not depend on the order of the matrix, the precomputation costs can be decreased considerably. This modification has been developed for trigonometric basis functions.

The costs of the computation step are minimal if the transformation matrix can be augmented by appending additional extra-elements, which is necessary to control the calculation accuracy. If the matrix elements are determined by orthogonal polynomials, to augment the matrix "on the right" and "from below” is not a problem. However, it is not always possible to append elements corresponding to negative values for the rows or columns, since in this case the functions that determine the values of the matrix elements may not be defined or change their behavior. If the matrix cannot be augmented, the block version of the extra-component method, which requires greater computational costs, should be used. The example of calculation of the Jacobi and Laguerre transforms has shown that the block version of the extra-component method compresses considerably the corresponding matrices. Multiplication by the compressed matrix requires additional FFT computations, but these constitute only a small fraction of the total computations. The above computational experiments have shown that the extra-component method can significantly decrease the computation time of matrix-vector multiplication for various classes of functions. The applicability of the approach proposed in this paper is not limited to these examples, and the developed algorithms can be adapted for other special matrices.
}

\section*{Acknowledgments}
The numerical implementation of proposed algorithms was carried out under state contract with ICMMG SB RAS (0251-2021-0004), the study of proposed algorithms  was financially supported by RFBR and Novosibirsk region  (Project No. 20-41-540003).

\section*{Conflict of Interest}
The authors declare that they have no conflict of interest.

\newpage
\normalsize
\bibliography{base}

\appendix
\setcounter{equation}{0}
\renewcommand\theequation{A.\arabic{equation}}

\clearpage\section{\textcolor{black}{Appendix A. An efficient precomputation procedure for trigonometric transforms}}
\color{black}
\normalsize

In the case of calculation of nonuniform trigonometric transformations, the precomputation costs of Algorithms~1 and~2 can be reduced from $O(N^2\log N)$ to $O(N \log N +Nq)$  arithmetic operations, where $q$  is determined by the required computation accuracy. As an example, let us consider an efficient procedure for calculating a row of the compressed matrix (\ref{A_comp}).

First, according to the inverse convolution theorem, for any $\mathbf{X},\mathbf{Y}\in \mathbb{C}^{N}$ we have
$$
\mathcal{F}\left(\mathbf{X} \odot \mathbf{Y} \right)= \mathcal{F}\mathbf{X} \ast \mathcal{F}\mathbf{Y}=\tilde{\mathbf{X}} \ast   \tilde{\mathbf{Y}},
$$
where $\odot$ denotes the component-wise product and
\begin{equation}
\left(\tilde{\mathbf{X}} \ast   \tilde{\mathbf{Y}}\right)_n\triangleq \frac{1}{\sqrt{N}}\sum_{m=0}^{N-1}\tilde{x}_m\tilde{y}_{(n-m)\text{mod}\ N}, \quad n=0,...,N-1.
\label{conv_appendix}
\end{equation}
Fig.~\ref{kaiser_win} shows that the Kaiser function can be approximated by a small number of coefficients of the Fourier series denoted here as $\tilde{x}_m$. Since in the discrete case for the Kaiser function no efficient formula for calculating individual Fourier components is known, the calculation of several components of the spectrum $\tilde{x}_m$  using FFT will require $O(N\log N)$ operations. This operation is performed once for all columns or rows of the transformation matrix.

Second, Fig.~\ref{pic:cheb_matrix22}b shows that there is no need to calculate all elements of a compressed matrix row by formula (\ref{conv_appendix}). It is sufficient to calculate about $24$ convolutions of the form (\ref{conv_appendix})  for each row to ensure an accuracy of transformation of the order $10^{-15}$.

Third, depending on the trigonometric transformation type given by the matrix $A$, the individual DFT components for the rows of the matrix $A$ can be calculated by one of the following formulas:
\begin{equation}
\begin{array}{ll}
\tilde{y}^{\mathrm{cos}}_k(\theta)=\frac{1}{\sqrt{N}}\sum_{j=0}^{N-1}\cos\left(j\theta\right)\omega_N^{jk} =\\\\ \frac{1}{\sqrt{N}}\frac{\omega_N^k\left((\left( \cos  \left(N \theta\right) -1)\cos \left(\theta\right)+\sin\left(\theta\right) \sin \left(N \theta\right)\right)-(\cos \left(N \theta\right)-1)\omega_N^{k}\right)}{\omega_N^{2k}-2 \omega_N^{k} \cos \left(\theta\right)+1},
\end{array}
\label{y1}
\end{equation}

\begin{equation}
\begin{array}{ll}
\tilde{y}^{\mathrm{sin}}_k(\theta)=\frac{1}{\sqrt{N}}\sum_{j=0}^{N-1}\sin\left(j\theta\right)\omega_N^{jk}\\\\=\frac{1}{\sqrt{N}}\frac{\omega_N^{k}\left(\left(\sin \left(\theta\right) (\cos \left(N \theta\right)-1)-\cos\left(\theta\right)  \sin\left(N\theta\right)\right)+\sin \left(N \theta\right)\omega_N^{k}\right)}{-\omega_N^{2k}+2 \omega_N^{k} \cos \left(\theta\right)-1},
\end{array}
\label{y2}
\end{equation}

\begin{equation}
\begin{array}{ll}
\tilde{y}^{\mathrm{exp}}_k(\theta)=\frac{1}{\sqrt{N}}\sum_{j=0}^{N-1}{\mathrm e}^{\mathrm{i} j\theta}\omega_N^{jk}=\frac{1}{\sqrt{N}}\frac{{\mathrm e}^{\mathrm{i} N \theta}-1}{{\mathrm e}^{\mathrm{i}\theta}\omega_N^{k}-1},
\end{array}
\label{y3}
\end{equation}
where $\omega_N=\exp{\left(-{2\pi\mathrm{i}}/{N}\right)}$. The variables $\theta$  and $k$ in formulas (\ref{y1})-(\ref{y3})  are separated, which reduces the precomputation time. As a result, the total cost of calculating all rows of the matrix  (\ref{A_comp})  will be proportional to  $O(N \log N+Nq)$, where the first term is the number of operations required for a single FFT calculation for the Kaiser function, and the second term is the number of operations of $N$ convolutions of the form (\ref{conv_appendix}) to calculate the compressed matrix elements.
\end{document}